\documentclass[12pt,leqno]{article}
\usepackage{amsmath}
\usepackage{amsthm}
\usepackage{amssymb}
\usepackage{xypic}

\input cyracc.def
        \newfam\cyrfam
        \font\tencyr=wncyr10 scaled 1200
        \font\sevencyr=wncyr10 scaled 700

        \textfont\cyrfam=\tencyr \scriptfont\cyrfam=\sevencyr
        \scriptscriptfont\cyrfam=\sevencyr

\makeatletter
\renewcommand{\subsubsection}{\@startsection
	{subsubsection}%
	{2}%
	{0mm}%
	{-\baselineskip}%
	{0.5\baselineskip}%
	{\normalfont\normalsize\itshape}}%
\makeatother
 
\makeatletter

\newcommand{\proofoftheoremname}{Proof of Theorem}
\newcommand{\remarkname}{Remark}
\makeatother

\newfam\cyrfam 
\font\tencyr=wncyr10 scaled 1200
\font\sevencyr=wncyr10

\textfont\cyrfam=\tencyr \scriptfont\cyrfam=\sevencyr
	\scriptscriptfont\cyrfam=\sevencyr

\def\CC{{{\bf C}}}
\def\QQ{{{\bf Q}}}
\def\ZZ{{{\bf Z}}}
\def\RR{{{\bf R}}}
\def\FF{{{\bf F}}}

\def\cD{{{\mathcal D}}} 
\def\cE{{{\mathcal E}}} 
\def\cF{{{\mathcal F}}} 

\def\cH{{{\mathcal H}}}

\def\cG{{{\mathcal G}}}
\def\cO{{{\mathcal O}}}

\def\bV{{{\mathbf{V}}}}
\def\bT{{{\mathbf{T}}}}

\def\bF{{{\mathbf{F}}}}

\def\cL{{{\mathcal L}}}
\def\cK{{{\mathcal K}}}

\def\cN{{{\mathcal N}}}

\def\cV{{{\mathcal V}}}
\def\cU{{{\mathcal U}}}

\newcommand{\hp}{H^1_{\text{par}}}

\def\epp{\varepsilon_{\gothp}}
\def\ep{\varepsilon}
\def\dpp{\delta_{\gothp}}
\def\DP{\Delta_{\gothp}}
\def\QSQ{\QQ_{\Sigma}/\QQ}
\def\LO{\Lambda_{\cO}}
\def\gothp{\mathfrak{p}}

\def\gotha{\mathfrak{a}}

\def\goths{\mathfrak{s}}

\def\gothm{\mathfrak{m}}

\def\gothq{\mathfrak{q}}

\def\Gal{\mathrm{Gal}}
\def\G{\mathrm{G}}
\def\Hom{\mathrm{Hom}}

\def\coker{\mathrm{coker}}
\def\Ker{\mathrm{Ker}}

\def\A{\mathrm{A}}

\def\Im{\mathrm{Im}}

\def\corank{\mathrm{corank}}
\def\Gal{\mathrm{Gal}}
\def\det{\mathrm{det}}

\def\rank{\mathrm{rank}}

\def\Ind{\mathrm{Ind}}

\def\ker{\mathrm{ker}}

\def\coker{\mathrm{coker}}
\def\dim{\mathrm{dim}}

\def\hm{\mathbf{h}_{\mathfrak m}}

\def\frob{\mathrm{Frob}}
\date{~}

\begin{document}

\title{Iwasawa theory for   Artin representations I} 
\author{Ralph Greenberg\thanks{Research supported in part by grants from the NSF.},  Vinayak Vatsal\thanks{Research supported in part by grants from the NSERC.} }
\date{}
\maketitle

\noindent {\em In honor of Professor Kenkichi Iwasawa}

\bigskip

\section{Introduction}

Suppose that $K$ is a finite Galois extension of $\QQ$.  Let 
$\Delta=\Gal(K/\QQ)$.  Consider an irreducible representation 
\[
\rho: ~ \Delta ~ \longrightarrow ~  GL_d(\overline{\QQ}) ~,
\]
where $\overline{\QQ}$ is an algebraic closure of $\QQ$ and $d \ge 1$.   Let $\chi$ denote the character of $\rho$.    If we fix an embedding 
$\sigma_{\infty}$  of $\overline{\QQ}$ into $\CC$,   then we obtain a $d$-dimensional representation of $\Delta$ over $\CC$. One can then define the corresponding Artin $L$-function 
$L_{\infty}(s, \chi)$.   It is defined by an Euler product for $Re(s)>1$  and can be analytically continued as a meromorphic function on the complex plane.  Conjecturally,  $L_{\infty}(s, \chi)$ is an entire function if $\chi$ is nontrivial. We will assume that $\chi$ is nontrivial throughout this paper. 

Now suppose that $p$ is a prime and that  we fix  an embedding $\sigma_p$  of   $\overline{\QQ}$ into $\overline{\QQ}_p$, an algebraic closure of the field $\QQ_p$ of $p$-adic numbers.  Then $\rho$ becomes a  $d$-dimensional representation of $\Delta$ over $\overline{\QQ}_p$. Its character $\chi$ now has values in   $\overline{\QQ}_p$.   One can ask whether there is a natural way to 
define a function $L_p(s, \chi)$  attached to $\rho$ and the fixed embedding,  where $s$ is  now a $p$-adic variable.    

This 
question has been studied previously when $K$ is totally real.  In that case, the values of $L_{\infty}(s, \chi)$ at negative odd integers are nonzero algebraic numbers.   Using the fixed embeddings $\sigma_{\infty}$ and $\sigma_p$, one can regard those values as elements of $\overline{\QQ}$ and then as elements of  $\overline{\QQ}_p$.  Such a $p$-adic $L$-function  $L_p(s, \chi)$ can then be defined by an interpolation property involving those elements of   $\overline{\QQ}_p$. First of all, if $d=1$,  then $\chi$ can be identified with an even Dirichlet character and  the construction of such a function $L_p(s, \chi)$  was done by Kubota and Leopoldt \cite{kl}.  
If $d > 1$, then there is a definition described in \cite{gre83} which is based on the Deligne-Ribet 
construction of $p$-adic $L$-functions for Hecke characters of the intermediate fields for the extension $K/\QQ$ given in \cite{dr}.    
Furthermore,  still assuming that $K$ is totally real,  one can give 
a precise interpretation of the zeros of  $L_p(s, \chi)$.  This is a 
generalization  of a famous conjecture of Iwasawa (for the case where $d=1$) and has been proven by Wiles \cite{wil90}.

Let $v$ denote an archimedian prime of $K$ and let $K_v$ denote the $v$-adic completion of $K$, which we identify with either $\RR$ or $\CC$.       As usual, we identify $\Gal(K_v/\RR)$ with a subgroup $\Delta_v$ of $\Delta$.   Let $d^{+}$ denote the multiplicity of the trivial character of $\Delta_v$ in $\rho|_{\Delta_v}$.  This does not depend on the choice of $v$.  Of course, if $K$ is totally real,  then $d^{+}=d$.  Our main objective
in this paper is to study   the case where $d^{+} = 1$ and $d \ge 1$.   Under certain assumptions, we will define certain $p$-adic $L$-functions and discuss a conjectural interpretation of the zeros of those functions.  Our theory reduces to the classical case where $d^{+}=d$ when $d=1$.   

We will usually make the following three-part assumption.   In the last part,  $\Delta_{\gothp}$ denotes the decomposition subgroup  of $\Delta$ for a prime $\gothp$ lying over $p$.  
\bigskip

\noindent {\bf Hypothesis A:} $~$  {\em The degree  $[K:\QQ]$ is not divisible by $p$. Also, $d^{+}=1$.   Furthermore, there exists  a 1-dimensional representation $\varepsilon_{\gothp}$ of $\Delta_{\gothp}$ which occurs with 
multiplicity 1 in  $\rho|_{\Delta_{\gothp}}$.} 
\bigskip

\noindent Our theory depends on choosing one such character $\varepsilon_{\gothp}$ of $\Delta_{\gothp}$.  Also, if $d > 1$, then $K$ must be totally complex and hence of even degree. We will therefore assume that $p$ is odd from here on.

Since $d^{+} = 1$, Frobenius Reciprocity implies that $\rho$ occurs with multiplicity 1 in the $\QQ$-representation of $\Delta$ induced from the trivial representation of $\Delta_v$ for any archimedian prime $v$ of $K$.  It follows that $\rho$ itself can be realized by matrices with entries in $\QQ(\chi)$, the field generated by the values of $\chi$ .  Let $\cF$ be the field generated over $\QQ_p$ by the values of $\chi$ and of $\varepsilon_{\gothp}$.   
We can then regard  $\rho$ as a representation of $\Delta$ over $\cF$  (via $\sigma_p$).   Let $V$ denote the underlying representation space. Thus,  $\dim_{\cF}(V)=d$.  Also,  the last part of Hypothesis A means that the maximal $\cF$-subspace  
$V^{(\varepsilon_{\gothp})}$  of $V$ on which $\Delta_{\gothp}$  acts by the character    $\varepsilon_{\gothp}$ has dimension 1.

The theory described in this paper depends not just on $\chi$, but also on the choice of the character $\varepsilon_{\gothp}$ of $\Delta_{\gothp}$ satisfying Hypothesis A.  Suppose that $\gothp'$ is another prime of $K$ lying over $p$.   
 Then $\gothp' = \delta(\gothp)$ for some $\delta \in \Delta$.  Thus,  $\Delta_{\gothp} = \delta^{-1} \Delta_{\gothp'}\delta$.      
Conjugation by $\delta^{-1}$ defines an 
isomorphism  $\Delta_{\gothp'} \to \Delta_{\gothp}$.  One then obtains a character   $\varepsilon_{\gothp'}$ of  $\Delta_{\gothp'}$  from   $\varepsilon_{\gothp}$ by composing with this isomorphism. One sees easily that $\varepsilon_{\gothp'}$ is well-defined and occurs with multiplicity 1 in $\rho|_{\Delta_{\gothp'}}$.  
In effect, we have a family of characters $\{\varepsilon_{\gothp}\}_{\gothp | p}$ for the decomposition subgroups which are compatible under conjugation.

Let $\cK$ denote the completion of $K$ at $\gothp$.  Then $\cK$ is a certain finite Galois  extension of $\QQ_p$. Consider the restriction map
   $\Gal(\cK/\QQ_p) \to \Delta_{\gothp}$. It is an isomorphism. 
Composing this map with the character $\varepsilon_{\gothp}$ of $\Delta_{\gothp}$  determines a 
character $\varepsilon$ of $\Gal(\cK/\QQ_p)$. It has values in $\cF$.  
If $\gothp'$ is another prime of $K$ lying over $p$,  then the completion $K_{\gothp'}$ is isomorphic to $\cK$.  The character $\varepsilon$ then determines a character  of $\Delta_{\gothp'}$.  In fact, this character does not depend on the isomorphism and is just $\varepsilon_{\gothp'}$ as defined above.   Thus, the family  $\{\varepsilon_{\gothp}\}_{\gothp | p}$ is determined by the single character  $\varepsilon$ of $\Gal(\cK/\QQ_p)$. 
\medskip

We will first describe the algebraic side of our theory. We assume that Hypothesis A is satisfied.  Let $\cO$ denote the ring of integers in $\cF$.  Let $T$ denote an $\cO$-lattice in $V$ which is invariant under the action of $\Delta$.  Since $|\Delta|$ is prime to $p$ and $\rho$ is irreducible over $\cF$,  $T$ is determined up to multiplication by an element of $\cF^{\times}$.  Let $D = V/T$. Thus, $D$ is isomorphic to 
$(\cF/\cO)^d$ as an $\cO$-module and has an $\cO$-linear action of $\Delta$. 
Furthermore,  for every prime  $\gothp$ of $K$ lying above $p$,  we will 
denote the image of $V^{(\varepsilon_{\gothp})}$  in $D$ by $D^{(\varepsilon_{\gothp})}$.   It is the maximal $\Delta_{\gothp}$-invariant  $\cO$-submodule of $D$  on which $\Delta_{\gothp}$ acts by $\varepsilon_{\gothp}$  and is isomorphic to  $\cF/\cO$ as an $\cO$-module.

We now define a {\em ``Selmer group''} associated to $\chi$ and $\varepsilon$ which we denote by  $S_{\chi, \varepsilon}(\QQ)$.  It will be a certain subgroup of 
$H^1(G_{\QQ}, D)$,  where  $G_{\QQ} = \Gal(\overline{\QQ}/\QQ)$.  As usual, we will write $H^1(\QQ, D)$ instead of $H^1(G_{\QQ}, D)$. A similar notation will be used for other fields and Galois modules. 
Suppose that  $\ell$ is a prime. Let  $\overline{\QQ}_{\ell}$ be an algebraic closure of $\QQ_{\ell}$.  Pick (arbitrarily)  an embedding of $\overline{\QQ}$ into $\overline{\QQ}_{\ell}$.  We can then identify $G_{\QQ_{\ell}} = \Gal(\overline{\QQ}_{\ell}/\QQ_{\ell})$  with a subgroup of 
$G_{\QQ}$.  Let $\QQ_{\ell}^{unr}$ be the maximal unramified extension of $\QQ_{\ell}$ in $\overline{\QQ}_{\ell}$.  The inertia subgroup of $G_{\QQ_{\ell}}$ is $G_{\QQ^{unr}_{\ell}}$ and is also identified with a subgroup of $G_{\QQ}$. 

Note that if $\ell=p$, then the chosen embedding of $\overline{\QQ}$ into $\overline{\QQ}_{p}$ determines a prime of $K$ lying above $p$, which we denote simply by $\gothp$.    Our Selmer group is defined by
\[
S_{\chi, \varepsilon}(\QQ) ~=~ \ker \Big( H^1(\QQ, D) ~\longrightarrow ~    \prod_{\ell \neq p} ~ H^1(\QQ_{\ell}^{unr}, D) \times H^1(\QQ_{p}^{unr}, D/D^{(\varepsilon_{\gothp})}) ~\Big) ~~.
\]
The {\em ``global-to-local''} map occurring in this definition is induced by the restrictions maps, where we identify the  $G_{\QQ^{unr}_{\ell}}$'s with subgroups of $G_{\QQ}$ as above.  However, 
as our notation suggests, $S_{\chi, \varepsilon}(\QQ)$ is completely determined by the $\cF$-valued character $\chi$ of $\Delta$ (of degree $d$) and by the $\cF$-valued character $\varepsilon$ of $\Gal(\cK/\QQ_p)$ (of degree 1).  It doesn't depend on the chosen embeddings and identifications. 

The above definition of our Selmer group is suggested by the general notion of a Selmer group which was studied in \cite{gre89} and which also depends on  having a filtration on a representation space   $V$ for $G_{\QQ}$ when restricted to the local Galois group $G_{\QQ_p}$.  In that paper, the local condition at the prime  $p$ is defined in terms of a certain subspace denoted by $F^{+}V$.  That paper primarily considers the so-called $p$-critical case in which the dimension of $F^{+}V$ is equal to  $d^{+}$.   And so, if $d^{+}=1$, then $F^{+}V$ should be of dimension 1.  In the definition given above, the role of $F^{+}V$ is now being played by $V^{(\varepsilon_{\gothp})}$.

Now $S_{\chi, \varepsilon}(\QQ)$ is a subgroup of $H^1(\QQ, D)$.  Since $D$ is an $\cO$-module, so is  $H^1(\QQ, D)$.  It is clear that  $S_{\chi, \varepsilon}(\QQ)$ is an $\cO$-submodule.  The restriction map 
\[
H^1(\QQ, D) ~\longrightarrow ~ H^1(K,  D)^{\Delta} ~=~ \Hom_{\Delta}( 
\Gal(K^{ab}/K),~D )
\]
is an isomorphism. Here $K^{ab}$ is the maximal abelian extension of $K$ in $\overline{\QQ}$.  Also, the notation $\Hom_{\Delta}(\cdot, ~\cdot)$ will always denote the group of continuous, $\Delta$-equivariant homomorphisms.  The above isomorphism follows from the assumption that $|\Delta|$ is not divisible by  $p$ and the fact that $G_K$ acts trivially on $D$.   Thus, we can identify $S_{\chi, \varepsilon}(\QQ)$ with its image under the restriction map which we will now describe. 

Let $\xi \in  \Hom_{\Delta}( \Gal(K^{ab}/K),~D)$.  
Let $K_{\xi}$ be the fixed field for $ker(\xi)$. Thus, $K_{\xi}$ is a finite, abelian extension of $K$, Galois over $\QQ$, and $\xi$ defines a  $\Delta$-equivariant isomorphism of   $\Gal(K_{\xi}/K)$ to a subgroup of $D$.  
For  $\xi$ to be in the image of $S_{\chi, \varepsilon}(\QQ)$, 
the local conditions defining $S_{\chi, \varepsilon}(\QQ)$ are equivalent to the following conditions on $\xi$.  First of all, 
the extension $K_{\xi}/K$ is unramified at all 
primes of $K$ not dividing $p$.  Furthermore, if $\gothp$ is a prime of $K$ lying over $p$,  then the image of the corresponding inertia subgroup of $\Gal(K_{\xi}/K)$ under the map $\xi$ is contained in $D^{(\varepsilon_{\gothp})}$.   
By using the above identification, we will prove the following result.  
\bigskip

\noindent {\bf Theorem 1.} {\em  $S_{\chi, \varepsilon}(\QQ)$ is finite.}
\bigskip

\noindent The proof of this theorem uses the Baker-Brumer theorem concerning the linear independence over $\overline{\QQ}$ of the $p$-adic logarithms of algebraic numbers and is reminiscent of the proof of Leopoldt's conjecture  for abelian extensions of $\QQ$.   Leopoldt's conjecture is usually formulated as a statement about the $p$-adic independence of units in a number field $K$.      But, if $K$ is Galois over $\QQ$, then  it  has an equivalent formulation as a statement about the 
{\em $\rho$-components} of the unit group, where $\rho$ varies over the absolutely irreducible Artin representations of $\Gal(K/\QQ)$.  As we will show in section 2, 
the Baker-Brumer theorem provides a proof of that conjecture  when $d^{+}=1$.  However, theorem 1 is a stronger statement and not just a consequence of Leopoldt's conjecture for such a $\rho$.  The proof requires a somewhat more refined application of the Baker-Brumer theorem and depends on Hypothesis A.

  Let $\QQ_{\infty}$ be the cyclotomic $\ZZ_p$-extension of $\QQ$.  For $n \ge 0$, let $\QQ_n$ denote the subfield of $\QQ_{\infty}$ of degree $p^n$ over $\QQ$.  One can define Selmer groups $S_{\chi, \varepsilon}(\QQ_{n})$  for all $n \ge 1$ and $S_{\chi, \varepsilon}(\QQ_{\infty})$ in a similar way to the definition of  $S_{\chi, \varepsilon}(\QQ)$  (the case $n=0$). We will actually prove that $S_{\chi, \varepsilon}(\QQ_{n})$ is finite for all $n \ge 0$.    Let $\Gamma = \Gal(\QQ_{\infty}/\QQ)$ and let $\Lambda_{\cO}=\cO[[\Gamma]]$.  Now $\Gamma$ acts naturally on $S_{\chi, \varepsilon}(\QQ_{\infty})$ and the action is $\cO$-linear. Thus, we can regard   $S_{\chi, \varepsilon}(\QQ_{\infty})$ as a discrete  $\Lambda_{\cO}$-module.  Let   $X_{\chi, \varepsilon}(\QQ_{\infty})$
 denote the Pontryagin dual of  $S_{\chi, \varepsilon}(\QQ_{\infty})$.   We can regard  $X_{\chi, \varepsilon}(\QQ_{\infty})$
 as a compact $\Lambda_{\cO}$-module.  The following result is  a straightforward consequence of the finiteness of the  $S_{\chi, \varepsilon}(\QQ_n)$'s. 
\bigskip

\noindent {\bf Theorem 2.} {\em The $\Lambda_{\cO}$-module $X_{\chi, \varepsilon}(\QQ_{\infty})$  is finitely-generated and torsion.  }
\bigskip

\noindent Alternatively, we would say that $S_{\chi, \varepsilon}(\QQ_{\infty})$ is a cofinitely-generated, cotorsion $\Lambda_{\cO}$-module.
\smallskip

      We will next state a result on the algebraic side concerning the characteristic ideal of  $X_{\chi, \varepsilon}(\QQ_{\infty})$. We will prove this result in a sequel to this paper (part II).  By definition, the characteristic ideal is a principal ideal in   $\Lambda_{\cO}$.   We will prove that it has a generator with an interpolation property of a certain form.  To simplify the statement and explanation, we will assume here that $\cF = \QQ_p$.  Thus, we are assuming that both $\chi$ and $\varepsilon$ have values in $\QQ_p$.  

Let $K_{\infty} = K\QQ_{\infty}$,  the cyclotomic $\ZZ_p$-extension of $K$.  We can identify $\Gal(K_{\infty}/K)$ with $\Gamma = \Gal(\QQ_{\infty}/\QQ)$ by the obvious restriction map,  noting that $K \cap \QQ_{\infty} = \QQ$ because $[K:\QQ]$ is prime to $p$.   We may assume  that $K$ contains $\mu_p$. Then $K_{\infty}=K(\mu_{p^{\infty}})$.  The action of  
$\Gal(K_{\infty}/K)$ on $\mu_{p^{\infty}}$ defines an isomorphism 
$\kappa: \Gamma \to 1+ p\ZZ_p$ (since $p$ is odd).   
Recall that $\cK$ denotes the completion of $K$ at any one of the primes of $K$ lying above $p$.   We let $\cK_{\infty}$ denote the cyclotomic $\ZZ_p$-extension of $\cK$.  We can also identify  $\Gal(\cK_{\infty}/\cK)$ with 
$\Gamma$.  For $n \ge 0$, the $n$-th layers in  $K_{\infty}/K$  and in $\cK_{\infty}/\cK$ will be denoted by $K_n$ and $\cK_n$, respectively.  We also will identify $\Gal(K_{\infty}/\QQ_{\infty})$ with $\Delta$ by the restriction map.  We then have a canonical isomorphism $\Gal(K_{\infty}/\QQ) \cong \Delta \times \Gamma$.  Similarly,  $\Gal(\cK_{\infty}/\QQ_p) \cong \cG \times \Gamma$, where $\cG = \Gal(\cK/\QQ_p)$.    We will denote $\Gal(K_n/\QQ)$ by $\Delta_n$. 

For each $n \ge 0$,  let $\cU_n$ denote the group of principal units in $\cK_n$. Let $U_n$ denote the group of units in $K_n$ which are principal in all the completions of $K_n$ at primes above $p$.  Thus, we can pick embeddings $K_n \to \cK_n$ in a compatible way and thereby obtain natural injective homomorphisms $U_n \to \cU_n$.  The lack of uniqueness will not affect the statement of the next theorem.    Those homomorphisms  commute with the natural global and local norm maps $U_m \to U_n$ and $\cU_m \to \cU_n$ for $m \ge n \ge 0$.   

 Since we are regarding $\cG$ as a subgroup of $\Gal(\cK_{\infty}/\QQ_p)$,  and the $\cU_n$'s are  abelian pro-$p$ groups and hence $\ZZ_p$-modules,   we can consider the $\cU_n$'s as  $\ZZ_p[\cG]$-modules.  For each $n \ge 0$, let $\cU_n^{(\varepsilon)}$ denote the $\varepsilon$-component  of $\cU_n$.  We can then define the projection map $\cU_n \to \cU_n^{(\varepsilon)}$, which is a surjective $\Gamma$-equivariant map.   We will use the somewhat peculiar notation $|\cdot|_{\varepsilon}$ for this projection map.    Our motivation for this notation is to suggest an analogy with the classical Stark conjecture which involves the complex log of the absolute value of so-called Stark units. Theorem 3 below involves instead the $p$-adic log of $|\cdot|_{\varepsilon}$ applied to certain units.

Suppose that $\varphi: \Gamma \to \overline{\QQ}_p^{\times}$ is a continuous group homomorphism.  Then we can extend $\varphi$ to a continuous $\ZZ_p$-algebra homomorphism $\Lambda \to \overline{\QQ}_p$,  which we also denote by $\varphi$.   Thus, its image $\varphi(\Lambda)$ is a compact subring of some finite extension of $\QQ_p$. The elements of $\Hom_{cont}(\Gamma,   \overline{\QQ}_p^{\times})$ of finite order will be of particular interest.  If $\varphi$ has order $p^n$,  then $\varphi(\Lambda) = \ZZ_p[\mu_{p^n}]$.    One can also consider the elements of $\Hom_{cont}(\Gamma,   \overline{\QQ}_p^{\times})$ of the form $\kappa^t$,  where $t \in \ZZ_p$.  They have values in $1+p\ZZ_p$ and we have $\kappa^t(\Lambda) = \ZZ_p$.

Now  $\varepsilon$ factors through the Galois group of a cyclic extension of $\QQ_p$ of degree prime to $p$. We have a canonical factorization $\varepsilon = \omega^a \beta$,  where $0 \le a < p-1$ and $\beta$ is unramified.  Let $b$ denote the order of $\beta$.  Thus, $\beta$ factors through $\Gal(\QQ_p(\mu_f)/\QQ_p)$ where $f=p^b-1$ and is a faithful character of that Galois group.   Also,  $\omega$ denotes the Teichm\"uller character and factors through $\Gal(\QQ_p(\mu_p)/\QQ_p)$.  If $\varphi$ is a character of $\Gamma$ of 
order $p^n$, then $\omega^a \varphi$ factors through $\Gal(\QQ_p(\mu_{p^{n+1}})/\QQ_p)$. Thus, $\ep\varphi$ can be regarded as a character of $\Gal(\QQ_p(\mu_{fp^{n+1}})/\QQ_p)$. We choose a generator $\zeta_f$ for $\mu_f$ so that its image $\tilde{\zeta}_f$ in the residue field $\FF_{p^b}$ for $\QQ_p(\mu_f)$ is part of a normal basis for $\FF_{p^b}/\FF_p$. Of course, this just means that $\{\tilde{\zeta}_f^{p^j}\}_{0 \le j<b}$ is linearly independent over $\FF_p$. The existence of such an element $\zeta_f$ is a theorem of Lenstra and Schoof (in \cite{LS}). For $n\ge 0$, we choose generators $\zeta_{fp^{n+1}}$ for $\mu_{fp^{n+1}}$ such that $\zeta_{fp^{n+1}}^p=\zeta_{fp^{n}}$.  We define a {\em ``Gaussian sum"} by 
\[
\goths(\ep\varphi) ~=~ \sum_{\sigma}~ \ep\varphi(\sigma^{-1}) \sigma(\zeta_{fp^{n+1}})
\]
where $\sigma$ varies over $\Gal(\QQ_p(\mu_{fp^{n+1}})/\QQ_p)$. The nonvanishing of $\goths(\ep\varphi)$ can be proved using our choice of $\zeta_f$.  The next theorem will not be proven here, but will be the main result in the sequel. 
It gives an interpolation property for a generator of the characteristic ideal. As we will discuss in the sequel, we can also give an interpolation property for $\varphi_{_0}(\theta_{\chi, \varepsilon})$ in most cases, where $\varphi_{_0}$ is the trivial character of $\Gamma$.

\bigskip

\noindent {\bf Theorem 3.}  {\em  Suppose that $d=2$ and that assumption A is satisfied. Then the characteristic ideal of $X_{\chi, \varepsilon}(\QQ_{\infty})$ has a generator $\theta_{\chi, \varepsilon}$ with the following property.  There exists a norm-compatible sequence $\{\eta_n\}$,  where  $\eta_n \in U_n$ for all $n \ge 0$,  such that   
\[
\varphi\big(\theta_{\chi, \varepsilon}\big) ~=~  \frac{1}{\goths(\ep\varphi)} \cdot \sum_{\delta \in \Delta_n} ~ \chi \varphi(\delta^{-1})\log_p\big( | \delta(\eta_n) |_{\varepsilon}\big) 
\]
for all  $n \ge 1$ and all $\varphi  \in  \Hom_{cont}(\Gamma,   \overline{\QQ}_p^{\times})$ of order $p^n$. }   
\bigskip

The proof of this theorem involves extending a classical result of Iwasawa to this setting. Indeed, we prove a certain four-term exact sequence involving the module  $X_{\chi, \varepsilon}(\QQ_{\infty})$ and certain inverse limits of local unit groups and global unit groups.   The sequence  $\{\eta_n\}$ defines an element $\eta_{\infty}$ in the $\Lambda$-module $\cU_{\infty} = \varprojlim~ \cU_n$.  The $\varepsilon$-component $\cU_{\infty}^{(\varepsilon)}$ has $\Lambda$-rank 1, and is usually just isomorphic to $\Lambda$.  The projection of $\eta_{\infty}$ to  $\cU_{\infty}^{(\varepsilon)}$ generates a $\Lambda$-submodule of  $\cU_{\infty}^{(\varepsilon)}$ of rank 1.  The corresponding quotient module has the same characteristic ideal as  $X_{\chi, \varepsilon}(\QQ_{\infty})$.  All of this will be discussed in the sequel. 
\medskip

We are not entirely satisfied with the above interpolation property. For one thing, it depends on the choice of $\zeta_f$. However, since we are free to multiply $\theta_{\chi,\ep}$ by an arbitrary unit in $\LO$,  we are hoping to find more natural interpolation factors, perhaps depending on $\chi$ and $\epsilon$ in some canonical way. In a previous version of this paper, we had the simpler (and canonical) factors $\tau(\omega^a \varphi)/p^{n+1}$, where $\tau(\omega^a \varphi)$ is the standard Gaussian sum.  However, that previous interpolation property corresponds to choosing a generator for the characteristic ideal in a large ring $\Lambda_{\cO'}$, where $\cO'$ is the ring of integers in a certain unramified extension $\cF'$ of $\cF$. 
\bigskip

Now we will discuss the analytic side of our theory.  Our most general results concern the case where $d=2$.    Since $d^{+}=1$,   we have $d^{-}=1$ and $\rho$ has odd determinant.  The action of $\Delta_{\gothp}$ on 
$V/V^{(\varepsilon_{\gothp})}$ is given by a character $\varepsilon'_{\gothp} :  
\Delta_{\gothp} \to \cF^{\times} $.   As before, we can define a character $\varepsilon':  \Gal(\cK/\QQ_p) \to \cF^{\times}$.  
 We assume the Artin conjecture for $\rho$ and all of its twists by 
 one-dimensional characters   so that we can associate to $\rho$ a certain modular form $f_{\rho}$ of weight 1.   We will define a $p$-adic $L$-function $L_p(s, \chi, \varepsilon)$ as the restriction of a 2-variable $p$-adic $L$-function to a certain line.  

For simplicity, we will assume here that $\varepsilon'$ is unramified. In general, we can reduce to that case by twisting by some  power of  $\omega$.  By using the fixed embeddings $\sigma_{\infty}$ and $\sigma_p$,  we can regard the $q$-expansion of $f_{\rho}$ as having coefficients in $\cF$.  One can then apply a theorem of Wiles to prove the existence of a Hida family of modular forms with the following property. They will be given by $q$-expansions with coefficients in  some finite extension of $\cF$.  Roughly speaking,  in each weight $k \ge 2$, the corresponding $p$-adic Galois representation has an unramified quotient (when considered as a representation of 
$G_{\QQ_p}$) and the action on that quotient  is given by a character $\varepsilon'_k$ of   $G_{\QQ_p}$.   
As the weight $k$ approaches 1 $p$-adically, $\varepsilon'_k$ approaches $\varepsilon'$. 

The two-variable $p$-adic $L$-function for a Hida family was first constructed by Mazur under certain assumptions and then by Kitagawa in \cite{kit94}, and later by a different method in \cite{gs}.  Under some mild restrictions on $\rho$, the Hida family is uniquely determined by the choice of $\varepsilon'$. However, since $\rho$ is 2-dimensional, the choice of $\varepsilon$ determines $\varepsilon'$, and we will use $\varepsilon$  in our notation. The corresponding two-variable $p$-adic $L$-function that we consider here is defined using the so-called canonical periods, following \cite{epw}.  The restriction to  $k=1$ defines a function $L_p(s,  \chi, \varepsilon)$ for $s \in \ZZ_p$.  Furthermore, there is an element $\theta^{an}_{\chi, \varepsilon}$ in $\Lambda_{\cO}$ such that 
\begin{equation}\label{eqn:LpDef}
L_p(s, \chi, \varepsilon) = \kappa^{1-s}(   \theta^{an}_{\chi, \varepsilon})
\end{equation}
for all $s \in \ZZ_p$. We conjecture that this element is a generator of the characteristic ideal of $X_{\chi, \varepsilon}(\QQ_{\infty})$.    In fact, when $d=1$, one of the equivalent forms of Iwasawa's classical Main Conjecture can be formulated in precisely this way.   (See the introduction of \cite{gre77}.)  In that case, 
one has $\epp = \chi |_{\Delta_{\gothp}}$ and $L_p(s, \chi, \ep)$ should be  taken to be $L_p(s, \chi)$, the Kubota-Leopoldt $p$-adic $L$-function for $\chi$.

We now denote the element $\theta_{\chi, \varepsilon}$ occurring in Theorem 3 by 
 $\theta^{al}_{\chi, \varepsilon}$.  Our conjecture asserts that  the principal ideals in $\Lambda_{\cO}$ generated by  
$\theta^{al}_{\chi, \varepsilon}$ and by  $\theta^{an}_{\chi, \varepsilon}$ are equal.   Neither  $\theta^{an}_{\chi, \varepsilon}$ nor  $\theta^{al}_{\chi, \varepsilon}$ is actually uniquely determined. They are only determined up to multiplication by an element of $\cO^{\times}$ for $\theta^{an}_{\chi, \varepsilon}$ or an element of $\Lambda_{\cO}^{\times}$ for $\theta^{al}_{\chi, \varepsilon}$.   In addition, although  $\theta^{al}_{\chi, \varepsilon}$ is nonzero, we cannot prove in general that  $\theta^{an}_{\chi, \varepsilon}$ is nonzero.

\bigskip

\section{Galois cohomology for Artin representations}

We will be quite general at first.  Suppose that $K$ is a finite Galois extension of $\QQ$, that $\Delta = \Gal(K/\QQ)$, and that  $\rho: \Delta \to GL_d(\cF)$  is any absolutely irreducible representation over $\cF$.   We even allow $[K:\QQ]$ to be divisible by $p$. Let $\chi$ denote the character of $\rho$. We use $\chi$ consistently in our notation.    We assume that $\cF$ is a finite extension of $\QQ_p$.  Let $\cO$ be the ring of integers in $\cF$.    Let $D = V/T$,  where $V$ is the underlying $\cF$-vector space for $\rho$ and $T$ is a $\Delta$-invariant $\cO$-lattice in $V$.    
 Let $\Sigma$ be any finite set of primes containing  $p$, $\infty$,  and the primes which are ramified in $K/\QQ$.  Thus, 
$K \subset \QQ_{\Sigma}$, where $\QQ_{\Sigma}$ is the maximal extension of $\QQ$ unramified outside of $\Sigma$.  We can regard $D$ as a discrete $\cO$-module with an $\cO$-linear action of $\Gal(\QQ_{\Sigma}/\QQ)$.  We denote $d$ by $d(\chi)$ and write $d(\chi) = d^{+}(\chi) + d^{-}(\chi)$, reflecting the action of the decomposition subgroup $\Delta_{v}$ for an archimedean prime $v$ of $K$. 
\medskip

We first discuss the Galois cohomology groups $H^i(\QQ_{\Sigma}/\QQ, D)$ for $i \ge 0$.  They are $\cO$-modules.  Of course, if $p \ge 3$,  then  $H^i(\QQ_{\Sigma}/\QQ, D)=0$ for $i \ge 3$.  For $p=2$, those cohomology groups are finite and have exponent $2$.  Thus, we just consider $i \in \{0,1,2\}$.  
We will be primarily  interested in the $\cO$-corank of  $H^1(\QQ_{\Sigma}/\QQ, D)$. The Poitou-Tate formula for the Euler-Poincar\'e characteristic gives
\[
 \corank_{\cO}\big(  H^1(\QQ_{\Sigma}/\QQ, D)\big) ~=~ d^{-}(\chi) + \corank_{\cO}\big(  H^2(\QQ_{\Sigma}/\QQ, D)\big) +  \corank_{\cO}\big(  H^0(\QQ_{\Sigma}/\QQ, D)\big)~~.
\]
If $\chi$ is the trivial character, then one sees easily that  $H^0(\QQ_{\Sigma}/\QQ, D)$ and  $H^1(\QQ_{\Sigma}/\QQ, D)$ both have $\cO$-corank  $1$, and hence  $H^2(\QQ_{\Sigma}/\QQ, D)$ has $\cO$-corank $0$.   If we assume that $\chi$ is not the trivial character,   then  $\corank_{\cO}\big(  H^0(\QQ_{\Sigma}/\QQ, D)\big)=0$.   The following result is an immediate consequence of the above remarks.
\bigskip

\noindent {\bf Proposition 2.1.}  {\em Assume that $\chi$ is nontrivial. Then $ \corank_{\cO}\big(  H^1(\QQ_{\Sigma}/\QQ, D)\big) \ge d^{-}(\chi)$ and equality holds if and only if  
$\corank_{\cO}\big(  H^2(\QQ_{\Sigma}/\QQ, D)\big) = 0$. }
\bigskip

In general, it is reasonable to conjecture that  $\corank_{\cO}\big(  H^2(\QQ_{\Sigma}/\QQ, D)\big) = 0$.  This is closely related to Leopoldt's conjecture, as we now explain.   Let
\[
\cU_p ~=~ \prod_{\gothp | p} ~ \cU_{\gothp} ~,
\] 
where $\cU_{\gothp}$ is the group of principal units in the completion $K_{\gothp}$.  Let  $U_K$ consists of the units of $K$  which are principal units at all $\gothp | p$.  Thus, $U_K$ is a subgroup of  the unit group of $K$ of finite index.
Consider the natural diagonal map 
\begin{equation}\label{eqn:diag}
U_K \longrightarrow   \cU_{p} ~~.
\end{equation}
   The $\cU_{\gothp}$'s, and hence $\cU_p$,  are $\ZZ_p$-modules. The map  (\ref{eqn:diag}) is obviously injective and extends to a map 
\begin{equation}\label{eqn:lambda-map}
\lambda_p: ~ U_K \otimes_{\ZZ}\ZZ_p \longrightarrow  \cU_p ~~.
\end{equation}
which is $\ZZ_p$-linear. The image of $\lambda_p$ is the closure of the image of  the map (\ref{eqn:diag}) and will sometimes be denoted simply by $\overline{U}_K$,   Furthermore, $\Delta$ acts naturally on both $U_K \otimes_{\ZZ} \ZZ_p$ and  $\cU_p$. The map $\lambda_p$ is $\Delta$-equivariant. Tensoring the above objects with $\cF$, we again get a $\Delta$-equivariant map which we denote by $\lambda_{p,\cF}$.   Then, restricting to the $\chi$-component,  we obtain a $\Delta$-equivariant map 
\begin{equation}\label{lambda-chi-map}
\lambda_{p,\cF}^{(\chi)}: ~  \big( U_K \otimes_{\ZZ} \cF \big)^{(\chi)}~ \longrightarrow~  \big( \cU_p \otimes_{\ZZ_p} \cF \big)^{(\chi)} ~~.
\end{equation}
Leopoldt's conjecture (for $K$ and $p$) asserts that $\lambda_p$ is injective. It is equivalent to saying that $\lambda_{p,\cF}$ is injective.   For each $\chi$, one can conjecture that $\lambda_{p,\cF}^{(\chi)}$ is injective.  This obviously follows from Leopoldt's conjecture.  Conversely, if  $\lambda_{p,\cF}^{(\chi)}$ is injective for all the absolutely irreducible character $\chi$ of $\Delta$, then Leopoldt's conjecture holds for $K$ and $p$.  Here one should choose $\cF$ so that 
all of the absolutely irreducible representations of $\Delta = \Gal(K/\QQ)$ are realizable over $\cF$. 
\medskip

Let $r(\chi)$ denote the multiplicity of $\chi$ in the $\Delta$-representation space  $U_K \otimes_{\ZZ} \cF$.   The proof of Dirichlet's Unit Theorem allows one to determine $r(\chi)$.  In fact, if $\chi$ is nontrivial, then $r(\chi) = d^{+}(\chi)$.   If $\chi$  is trivial, then $r(\chi)=0$.  On the other hand, the $\Delta$-representation space  $\cU_p \otimes_{\ZZ_p} \cF$ is isomorphic to the regular representation of $\Delta$ over $\cF$ and hence $\chi$ has multiplicity $d(\chi)$.   The image of $\lambda_{p,\cF}^{(\chi)}$ is the $\chi$-component of $\Im(\lambda_{p,\cF})$.  Let $r_{p}(\chi)$ denote the multiplicity of $\chi$ in  $\Im(\lambda_{p,\cF})$.  The conjecture that  $\lambda_{p,\cF}^{(\chi)}$ is injective can be equivalently stated as follows.
\bigskip

\noindent {\bf LC$(\chi,p)$:} $~~$   {\em We have $r_{p}(\chi)  =  r(\chi)$.}
\bigskip

\noindent It is obvious that $r_{p}(\chi)  \le r(\chi)$.  One can reformulate  {\bf LC$(\chi,p)$}   in terms of Galois cohomology. 
\bigskip 

\noindent {\bf Proposition 2.2.}  {\em  The assertion  {\bf LC$(\chi,p)$} holds if and only if  $\corank_{\cO}\big(  H^2(\QQ_{\Sigma}/\QQ, D)\big)=0$.}
\bigskip

\noindent {\em Proof.}    Consider the restriction map 
\[  
   \alpha:  ~ H^1(\QQ_{\Sigma}/\QQ, D)  \longrightarrow   H^1(\QQ_{\Sigma}/K, D)^{\Delta} ~~.
\]
Since $\Delta$ is finite, it follows that the kernel and cokernel of $\alpha$ are both finite.  Hence  $H^1(\QQ_{\Sigma}/\QQ, D)$ and   $H^1(\QQ_{\Sigma}/K, D)^{\Delta}$ have the same $\cO$-corank.  Furthermore, 
\[
H^1(\QQ_{\Sigma}/K,~ D)^{\Delta} ~=~ \Hom_{\Delta}( \Gal(M/K),~ D)~~,
\]
where $M$ is the maximal abelian pro-$p$ extension of $K$ contained in $\QQ_{\Sigma}$.  Clearly, $M$ is Galois over $\QQ$ and so $\Delta$ acts naturally on $\Gal(M/K)$.  Also, $\Gal(M/K)$ can be considered as a $\ZZ_p$-module and is 
finitely generated.   

Let $\widetilde{K}$ denote the compositum of all $\ZZ_p$-extensions of $K$.  Then $K \subset \widetilde{K} \subseteq M$ and $M/\widetilde{K}$ is a finite extension.  Furthermore, $\widetilde{K}$ is Galois over $\QQ$ and hence $\Delta$ also acts naturally on $\Gal(\widetilde{K}/K)$.  These remarks imply the first of the following equalities.  The second equality is easily verified.   
\[  \corank_{\cO}\big(H^1(\QQ_{\Sigma}/\QQ, ~D)\big)   ~=~  \corank_{\cO}\big(\Hom_{\Delta}( \Gal(\widetilde{K}/K)\otimes_{\ZZ_p}\cO , ~D) \big)   
\]
\[
~=~ \dim_{\cF}\big( \Hom_{\Delta}(  \Gal(\widetilde{K}/K)\otimes_{\ZZ_p}\cF , ~V)\big) ~~. 
\]
Hence, the $\cO$-corank of $H^1(\QQ_{\Sigma}/\QQ, D)$  is just   
the multiplicity of $\chi$ in the $\Delta$-representation space  $\Gal(\widetilde{K}/K)\otimes_{\ZZ_p}\cF$.   We next describe that multiplicity in terms of 
 $r_{p}(\chi)$.  

By class field theory, there is a homomorphism 
\[
 \cU_p \big/\Im(\lambda_p)~ \longrightarrow ~ \Gal(\widetilde{K}/K)
\]
whose kernel and cokernel are finite.  The multiplicity of $\chi$ in  $ \cU_p \otimes_{\ZZ_p} \cF$ is $d(\chi)$.   The multiplicity of $\chi$  in $\Im(\lambda_p)\otimes_{\ZZ_p}\cF$ is  $r_{p}(\chi)$.  Hence the multiplicity of $\chi$ in $ \Gal(\widetilde{K}/K)\otimes_{\ZZ_p}\cF$ is $d(\chi) -  r_{p}(\chi)$.  

On the other hand, as mentioned above,  if $\chi \neq \chi_0$, then one has $r(\chi) = d^{+}(\chi)$.  
This follows from the well-known fact that 
\[
U_K \otimes_{\ZZ} \QQ ~\cong ~ \Ind_{\Delta_v}^{\Delta} (\theta_0)\big/ V_{0}
\]
where $v$ is an archimedean prime of $K$,  $\theta_0$ is the trivial character of the decomposition subgroup $\Delta_v$, and $V_{0}$ is the underlying $\QQ$-representation space for the trivial character $\chi_0$ of $\Delta$ over $\QQ$.  If $\chi \neq \chi_0$, then Frobenius Reciprocity implies that   $r(\chi)$ is equal to the multiplicity of $\theta_0$ in $\rho |_{\Delta_v}$, which is indeed  equal to $d^{+}(\chi)$.   Of course,  we have $r(\chi_0) = 0$.

 Using the formula for the  Euler-Poincar\'e characteristic, we obtain
\[ 
  \corank_{\cO}\big(H^2(\QQ_{\Sigma}/\QQ, ~D)\big) ~=~  \big(d(\chi) -  r_{p}(\chi) 
\big) -   d^{-}(\chi) -   \corank_{\cO}\big(H^0(\QQ_{\Sigma}/\QQ, ~D)\big) ~~.
\]
If $\chi \neq \chi_0$, then  we obtain $ \corank_{\cO}\big(H^2(\QQ_{\Sigma}/\QQ, ~D)\big) = r(\chi) - r_{p}(\chi)$.  This gives the stated equivalence in proposition 2.1.  For $\chi = \chi_0$,   the result is obvious.     \hfill$\blacksquare$
\bigskip

\noindent {\bf Remark 2.3.}  If $p$ is an odd prime, then the $p$-cohomological dimension of $\Gal(\QQ_{\Sigma}/\QQ)$ is 2. Hence $H^3(\QQ_{\Sigma}/\QQ, ~D[p])=0$.
Since $D$ is divisible by $p$, it follows that $H^2(\QQ_{\Sigma}/\QQ, ~D)$ is also divisible by $p$.  Thus, its Pontryagin dual is a finitely-generated, torsion-free $\ZZ_p$-module, and hence is free. It is free as an $\cO$-module. In particular, if the $\cO$-corank of  $H^2(\QQ_{\Sigma}/\QQ, ~D)$ is 0, then we actually have  $H^2(\QQ_{\Sigma}/\QQ, ~D)=0$.   If $p=2$, one can just show that $pH^2(\QQ_{\Sigma}/\QQ, ~D)$ is 0.
 \bigskip
 
One useful positive result concerning  {\bf LC$(\chi,p)$} is the following. 
\bigskip

\noindent {\bf Proposition 2.4.}  {\em  If   $r(\chi) \ge 1$, then $r_{p}(\chi) \ge 1$.}
\bigskip

\noindent The proof makes use of the following lemma.   We regard 
$U_K\otimes_{\ZZ}\overline{\QQ}$ as a module over the group ring $\overline{\QQ}[\Delta]$, where $\overline{\QQ}$ is the algebraic closure of $\QQ$ in $\overline{\QQ_p}$.  We regard $\Im(\lambda_p)\otimes_{\ZZ_p}\overline{\QQ}_p$ as a module over $\overline{\QQ}_p[\Delta]$ and hence over the subring  $\overline{\QQ}[\Delta]$.
\bigskip

\noindent {\bf Lemma 2.5.}  {\em  Suppose that 
$\theta \in \overline{\QQ}[\Delta]$ and that $\theta$ does not annihilate $U_K\otimes_{\ZZ}\overline{\QQ}$.   Then $\theta$ does not annihilate  $\Im(\lambda_p)\otimes_{\ZZ_p}\overline{\QQ}_p$.}
\bigskip

\noindent {\em Proof.}   Let $\cF$ be a finite extension of $\QQ_p$ containing the coefficients of $\theta$.   As before, we let $\cK$ be the completion of $K$ at one of the primes of $K$ lying above $p$.  Let $\cU$ be the group of principal units in $\cK$.     We consider the corresponding projection map 
\[
pr:~     \cU_p \otimes_{\ZZ_p}\cF ~= ~\prod_{\gothp | p}  \big(\cU_{\gothp} \otimes_{\ZZ_p}\cF\big) ~ \longrightarrow ~  \cU \otimes_{\ZZ_p}\cF 
\]
It suffices to  show that $\Im(\lambda_{p,\cF})^{\theta}$ has a nontrivial image under the map 
$pr$. 
\medskip

By assumption, there exists an $\eta \in U_K$ such that $\eta^{\theta}$ is nontrivial (as an element in  $U_K\otimes_{\ZZ}\cF$).   Now 
$\theta = \sum_{\delta \in \Delta} \alpha_{\delta}\delta$, where the $\alpha_{\delta}$'s are in $\cF \cap \overline{\QQ}$.  Choose a maximal subset $\Delta'$ of $\Delta$ so that $\{\delta'(\eta)\}_{\delta' \in \Delta'}$ is a multiplicatively independent set.  Thus, replacing $\theta$ by an integral multiple if necessary, we have 
\[
\eta^{\theta} ~=~ \eta^{\psi}
\]
where $\psi = \sum_{\delta' \in \Delta'} \beta_{\delta'} \delta'$ and the   $\beta_{\delta}$'s are in $\cF \cap \overline{\QQ}$. Since $\eta^{\psi}$ is nontrivial, the  $\beta_{\delta}$'s are not all zero.  
\medskip

 The map  $log_p: \cU \to \cK$ extends to an $\cF$-linear map from  
$\cU \otimes_{\ZZ_p}\cF$ to   $\cK \otimes_{\ZZ_p}\cF$.   The image of $\eta^{\psi}$ under the map 
\[
log_p \circ pr \circ \lambda_{p,\cF} : ~   U_K\otimes_{\ZZ}\cF ~\longrightarrow ~  \cK \otimes_{\ZZ_p}\cF~~.
\]
is  
$\gamma = \sum_{\delta' \in \Delta'} log_p\big(\delta'(\eta)\big)\otimes \beta_{\delta'}$.  The inclusions $\cK, \cF \subset \overline{\QQ}_p$ induce a map from  $\cK \otimes_{\ZZ_p}\cF$ to  $\overline{\QQ}_p$.  The image of $\gamma$ is $ \sum_{\delta' \in \Delta'} \beta_{\delta'} log_p\big(\delta'(\eta)\big)$.  The Baker-Brumer theorem \cite{brumer}   implies that $\gamma \neq 0$.  It follows that the image of $\eta^{\theta}$ under $\lambda_{p,\cF}$ is nontrivial. Since $\lambda_{p,\cF}$ is $\Delta$-equivariant, it indeed follows that $\theta$ is not an annihilator of   $\Im(\lambda_p)\otimes_{\ZZ_p}\overline{\QQ}_p$.   \hfill$\blacksquare$
 \bigskip

\noindent {\bf Remark 2.6.}   The proof of the above lemma shows that if $\theta \in \cF[\Delta]$ has coefficients in $\cF\cap\overline{\QQ}$ and does not annihilate $U_K\otimes_{\ZZ}\cF$, then $pr\big( \Im(\lambda_{p,\cF})^{\theta}\big)$ is a nontrivial $\cF$-subspace of 
$\cU\otimes_{\ZZ_p}\cF$.   This will be important later. 
\bigskip

Proposition 2.4 follows from lemma 2.5.  One takes $\theta$ to be the 
idempotent for $\chi$ in $\cF[\Delta]$ which does have coefficients in $\cF \cap \overline{\QQ}$.  The following corollary follows immediately since $r(\chi) =  d^{+}(\chi)$ when $\chi \neq \chi_0$. 
\bigskip

\noindent {\bf Corollary 2.7.}  {\em  If $d^{+}(\chi)=1$,  then {\bf LC$(\chi,p)$} is true.}
\bigskip

   A very closely related result is proved in \cite{ekw}. The above corollary is implicit in that work. In addition, those authors give an interesting class of examples where Leopoldt's conjecture can be proved. In fact,  their examples are an illustration of the above corollary. Suppose that $\Delta \cong A_4$, the alternating group of order 12. Suppose also that $K$ is not totally real. Then $U_K$ has rank 5. There are four absolutely irreducible representations of $\Delta$, up to isomorphism, three of dimension 1 and one of dimension 3. One finds easily that the three nontrivial representations are the constituents in $U_K\otimes_{\ZZ}\cF$ as a representation space for $\Delta$, two of dimension 1 and the one of dimension 3, and that they all satisfy $d^{+}=1$. And so Leopoldt's conjecture for $K$ and $p$ follows from the above corollary.
\bigskip

 The rest of this section concerns the Galois cohomology groups 
 $H^i(\QQ_{\Sigma}/\QQ_{\infty}, D)$.  They are discrete $\cO$-modules with a continuous action of $\Gamma = \Gal(\QQ_{\infty}/\QQ)$ and so can be regarded as discrete $\Lambda_{\cO}$-modules,  where $\Lambda_{\cO} = \cO[[\Gamma]]$.  They are cofinitely generated as $\Lambda_{\cO}$-modules. Their $\Lambda_{\cO}$-coranks can be determined without assuming any conjectures. The results below are essentially consequences of theorems of Iwasawa, notably theorems 17 and 18 in \cite{iwa73}.  We also refer to \cite{gre89} and \cite{gre06}. 
  \bigskip

\noindent {\bf Proposition 2.8.}  {\em The $\Lambda_{\cO}$-corank of   $H^1(\QQ_{\Sigma}/\QQ_{\infty}, D)$ is equal to $d^{-}(\chi)$.  If $p$ is odd, then  $H^2(\QQ_{\Sigma}/\QQ_{\infty}, D) = 0$. If $p=2$, then  $pH^2(\QQ_{\Sigma}/\QQ_{\infty}, D)=0$. } 
\bigskip

\noindent {\em Proof.}    The theorems of Poitou and Tate imply a formula for the  Euler-Poincar\'e characteristic involving the $\Lambda_{\cO}$-coranks of  the $H^i(\QQ_{\Sigma}/\QQ_{\infty}, D)$'s. The alternating sum of those $\Lambda_{\cO}$-coranks for $0 \le i \le 2$  is $-d^{-}(\chi)$.  (See proposition 3 in \cite{gre89}.)    Since the $\Lambda_{\cO}$-corank of $H^0(\QQ_{\Sigma}/\QQ_{\infty}, D)$ is clearly $0$,  one obtains the formula
\[
 \corank_{\Lambda_{\cO}}\big(  H^1(\QQ_{\Sigma}/\QQ_{\infty}, D)\big) ~=~ d^{-}(\chi) + \corank_{\Lambda_{\cO}}\big(  H^2(\QQ_{\Sigma}/\QQ_{\infty}, D)\big)~~.
\]
 The proposition amounts to the assertion that  $H^2(\QQ_{\Sigma}/\QQ_{\infty}, D)$ has $\Lambda_{\cO}$-corank 0.  
\medskip

Let $K_{\infty}=K\QQ_{\infty}$, the cyclotomic $\ZZ_p$-extension of $K$.  Consider the restriction map 
\[
 H^2(\QQ_{\Sigma}/\QQ_{\infty},~ D) \longrightarrow H^2(\QQ_{\Sigma}/K_{\infty}, ~D) ~~.
\]
Note that $\Gal(\QQ_{\Sigma}/K_{\infty})$ acts trivially on $D$.  Thus, 
as a  $\Gal(\QQ_{\Sigma}/K_{\infty})$-module, $D$ is isomorphic to $D_{\chi_0}^{d^{-}(\chi)}$. 
Furthermore, the kernel of the above map has exponent dividing the degree $[K_{\infty}:\QQ_{\infty}]$ and hence has $\Lambda$-corank 0.    

Now $H^1(\QQ_{\Sigma}/K_{\infty}, ~D_{\chi_0}) = \Hom( \Gal(M^{\Sigma}_{\infty}/K_{\infty}), D_{\chi_0})$.  Here $M^{\Sigma}_{\infty}$ denotes the maximal, abelian pro-$p$ extension of $K_{\infty}$ contained in $\QQ_{\Sigma}$. One can regard $\Gal(M^{\Sigma}_{\infty}/K_{\infty})$ as a $\Lambda$-module.  
According to a well-known theorem of Iwasawa, the $\Lambda$-rank of  $\Gal(M^{\Sigma}_{\infty}/K_{\infty})$ is equal to $r_2(K)$, the number of complex primes of $K$.   (See Theorem 17 in \cite{iwa73} for the case where $\Sigma$ is the set of primes $\Sigma_0$ lying above $p$ or $\infty$. For  larger $\Sigma$,   one can use the fact that there are only finitely many primes of $K_{\infty}$ lying above primes in $\Sigma-\Sigma_0$ 
to show that $\Gal(M^{\Sigma}/M^{\Sigma_0})$ has finite $\ZZ_p$-rank and hence is a torsion $\Lambda$-module.  Thus, the $\Lambda_{\cO}$-corank of  $H^1(\QQ_{\Sigma}/K_{\infty}, ~D_{\chi_0})$ is $r_2(K)$.  Again, the Poitou-Tate Duality Theorems imply that the Euler-Poincar\'e characteristic (which involves the alternating sum of the $\Lambda_{\cO}$-coranks of the  $H^i(\QQ_{\Sigma}/K_{\infty}, D)$'s for $0 \le i \le 2$)  is equal to $-r_2(K)$.  It follows that $H^2(\QQ_{\Sigma}/K_{\infty}, ~D_{\chi_0})$  has $\Lambda_{\cO}$-corank 0. Therefore,  $H^2(\QQ_{\Sigma}/\QQ_{\infty}, D)$ indeed has $\Lambda_{\cO}$-corank 0.

If $p$ is odd, then proposition 4 in \cite{gre89} asserts that the Pontryagin dual of $H^2(\QQ_{\Sigma}/\QQ_{\infty}, D)$ is a free $\Lambda$-module. Since its rank is 0, it must vanish. Hence  $H^2(\QQ_{\Sigma}/\QQ_{\infty}, D)=0$.     By considering the restriction map to an imaginary quadratic extension of $\QQ_{\infty}$, the argument in \cite{gre89} shows that $H^2(\QQ_{\Sigma}/\QQ_{\infty}, D)$ has exponent 1 or 2 when $p=2$.     \hfill$\blacksquare$
\bigskip

One other basic result about the structure of the $\Lambda_{\cO}$-module $H^1(\QQ_{\Sigma}/\QQ_{\infty}, D)$ concerns the phenomenon of {\em purity}, namely that the torsion submodule of the Pontryagin dual of that $\Lambda_{\cO}$-module has support which is purely of codimension 1.  This is the content of the following proposition. 
\bigskip

\noindent {\bf Proposition 2.9.}  {\em Assume that $p$ is odd. Then the Pontryagin dual of $H^1(\QQ_{\Sigma}/\QQ_{\infty}, D)$ has no nontrivial, finite $\Lambda_{\cO}$-submodules. }
\bigskip

\noindent {\em Proof.}  This follows immediately from proposition 5 in \cite{gre89}, the crucial assumption being that $H^2(\QQ_{\Sigma}/\QQ_{\infty}, D)$ vanishes.    If the image of $\rho$ has order prime to $p$, one could alternatively use the fact that the $\Lambda$-module $\Gal(M^{\Sigma}_{\infty}/K_{\infty})$ has no nontrivial, finite $\Lambda$-submodules. This is a consequence of theorem 18 in \cite{iwa73}.    \hfill$\blacksquare$

\bigskip

\bigskip

\section{ Selmer groups for Artin representations 1}  We will prove theorem 1 in this section. More generally, we prove the finiteness of the Selmer group over $\QQ_n$ for all $n \ge 0$. 
Suppose that $\rho$ is an  irreducible Artin representation of $G_{\QQ}$ over $\cF$ satisfying Hypothesis A in the introduction.   Since $d^{+}=1$, one sees easily that $\rho$ is absolutely irreducible.   We assume that $\rho$ factors through $\Delta=\Gal(K/\QQ)$, where $K$ is a finite Galois extension of $\QQ$  of degree prime to $p$.  We consider the Selmer groups $S_{\chi, \varepsilon}(\QQ_{n})$ for all $n \ge 0$.   The definition of the above Selmer groups can be described as follows.  As usual, the description of the local condition at a prime $\nu$ of a number field
of residue field characteristic $\ell$  always implicitly involves choosing a fixed embedding of $\overline{\QQ}$ into $\overline{\QQ}_{\ell}$ which induces the prime $\nu$ on the number field.  The Selmer group is independent of those choices. 
\medskip

If $\ell$ is any finite prime and $\nu$ is a prime of $\QQ_n$ lying over $\ell$, we denote the completion of $\QQ_n$ at $\nu$ by $\QQ_{n, \nu}$. If $\nu$ is a prime of $\QQ_{\infty}$ lying over $\ell$, then we let $\QQ_{\infty, \nu} =  \bigcup_{n \ge 0} \QQ_{n, \nu}$. Note that if $\ell \neq p$, then  $\QQ_{\infty, \nu} \subset \QQ_{\ell}^{unr}$.  For such a prime $\ell$, we define 
\[
\cH_{\ell}(\QQ_{n}, D) = \prod_{\nu | \ell} ~ \Im \left( H^1(\QQ_{n, \nu}, D) \longrightarrow  H^1(\QQ_{\ell}^{unr},  D)\right)~~
\]
for all $n \ge 0$.  There is a unique prime $\pi$ of $\QQ_{\infty}$  lying over $p$. Thus, our chosen embedding of $\overline{\QQ}$ into $\overline{\QQ}_{p}$ induces $\pi$ on $\QQ_{\infty}$. Let $\gothp$ be the prime of $K$ induced by that embedding. We then define 
\begin{equation}\label{eqn:cH_p}
 \cH_{p}(\QQ_{n}, D)  = \Im \left(   H^1(\QQ_{n, \pi},~ D) \longrightarrow H^1(\QQ_{n, \pi}^{unr},  ~ D/D^{(\varepsilon_{\gothp})}\right)~.
\end{equation}
We can then define the Selmer groups to be considered by
\begin{equation}\label{eqn:SelDefn}
S_{\chi, \varepsilon}(\QQ_{n}) =   \Ker \Big( H^1(\QQ_{\Sigma}/\QQ_{n}, D) ~\longrightarrow ~\prod_{\ell \in \Sigma} ~ \cH_{\ell}(\QQ_{n}, D) \Big)~.
\end{equation}
for any $n \ge 0$. Here $\Sigma$ is a finite set of primes of $\QQ$ containing $p$  and $\infty$ chosen so that  $K \subset \QQ_{\Sigma}$. The global-to-local map occurring in this definition is defined by the various restriction maps.  Since we are taking $p$ to be odd, the 1-cocycle classes are certainly trivial at the
archimedean primes of $\QQ_n$.  And so the product in (\ref{eqn:SelDefn}) is just over the nonarchimedean primes in $\Sigma$.   

Note that if $\chi$ is trivial, then $\epp$ is trivial and $D^{(\epp)} = D$. It is then obvious that 
\begin{equation}\label{eqn:trivial}
\Hom(\Gal(\QQ_{\infty}/\QQ_{n}), D)~ \subseteq ~S_{\chi, \ep}(\QQ_n)
\end{equation}
 and this Selmer group is therefore infinite. Actually, one can easily show that equality holds in (\ref{eqn:trivial}). However, as we stated in the introduction, we are assuming throughout this paper that $\chi$ is nontrivial.   We will now prove the following result.  
\bigskip

\noindent {\bf Proposition 3.1.} {\em Assume that Hypothesis A is satisfied. Then $S_{\chi, \varepsilon}(\QQ_{n})$ is finite for all $n \ge 0$.}
\bigskip

\noindent {\em Proof.}   Since $p \nmid [K:\QQ]$, we have $K \cap \QQ_n = \QQ$.   Let $K_n = K\QQ_n$, the $n$-th layer of the cyclotomic $\ZZ_p$-extension of $K$.       Then we can identify $\Gal(K_n/\QQ_n)$ with  $\Delta$ and $\Gal(K_n/K)$ with $\Gal(\QQ_n/\QQ)$, which in turn can be identified with $\Gamma/\Gamma_{n}$. Here 
$\Gamma_n=\Gamma^{p^n}=\Gal(\QQ_{\infty}/\QQ_n)$.  
 Let $\Delta_n=\Gal(K_n/\QQ)$.  We have a canonical isomorphism $\Delta_n \cong \Delta \times \Gamma/\Gamma_n$.    Each prime $\gothp$ of $K$ lying over $p$ is totally ramified in $K_n/K$.  Let $\gothp_n$ be the unique prime of $K_n$ lying over $\gothp$ and let $\cU_{n,\gothp}$ denote the group of principal units in the completion of $K_n$ at $\gothp_n$.     The decomposition subgroups of $\Gal(K_n/\QQ_n)$ and $\Delta_n$ for $\gothp_n$ are then identified   with $\Delta_{\gothp}$ and $\Delta_{\gothp} \times \Gamma/\Gamma_n$, respectively.  Both of these decomposition groups act on $\cU_{n,\gothp}$.  There is a natural action of $\Delta_n$ on 
 \[
 \cU_{n,p} ~=~ \prod_{\gothp | p} ~ \cU_{n, \gothp}~~.
 \]
As in section 2, we have a natural $\Delta_n$-equivariant map 
\[
\lambda_{n,p}: ~ U_{n}\otimes_{\ZZ}\ZZ_p~\longrightarrow~ \cU_{n,p}~~,
\]
where $U_{n}$ is the group of units of $K_n$ which are principal units at $\gothp_n$ for all $\gothp | p$. The image of $\lambda_{n,p}$ will be denoted by $\overline{U}_{n,p}$, 
\smallskip

We consider the Selmer groups $S_{\chi, \varepsilon}(\QQ_{n})$      for $n \ge 0$. 
The restriction map 
\[
H^1(\QQ_n, D) ~  \longrightarrow ~ H^1(K_n, D)^{\Delta} ~=~ \Hom_{\Delta}(G_{K_n}^{ab}, D)~~,
\]
is an isomorphism and the image of $S_{\chi, \varepsilon}(\QQ_{n})$ is contained in $\Hom_{\Delta}(\Gal(M_n/K_n), D)$,  where $M_n$ is the maximal abelian pro-$p$ extension of $K_n$ unramified outside of the set  of primes lying over $p$.  This assertion captures the local conditions defining $S_{\chi, \varepsilon}(\QQ_{n})$ at all primes $\ell \neq p$. 
Suppose that $\xi \in  \Hom_{\Delta}(\Gal(M_n/K_n), D)$. Noting that the local degree 
$[K_{n, \gothp_n} : \QQ_{n,\pi}]$ is prime to $p$, one sees that
 $\xi$ is in the image of 
$S_{\chi, \varepsilon}(\QQ_{n})$ if and only if  the image of the inertia subgroup $I_{n, \gothp}$ of $\Gal(M_n/K_n)$ for $\gothp_n$ under $\xi$ is contained in  $D^{(\varepsilon_{\gothp})}$. Note that if this condition is satisfied for one $\gothp |p$, it is satisfied for all those $\gothp$'s.   Also, $I_{n, \gothp}$ is 
$\Delta_{\gothp}$-invariant.  Now we can extend the map $\xi$ to a $\Delta$-equivariant map
\[
\xi_{\cO}:   \Gal(M_n/K_n)_{\cO} ~ \longrightarrow ~ D~~,
\]
where the subscript $\cO$ here (and elsewhere) denotes the tensor product over $\ZZ_p$ with $\cO$, a ring which is finite and flat over $\ZZ_p$.  The map $\xi_{\cO}$ is continuous and $\cO$-linear.  
Furthermore, with this notation, we have a $\Delta_{\gothp}$-decomposition
\[
 I_{n, \gothp, \cO}  ~=~   I_{n,\gothp, \cO}^{(\varepsilon_{\gothp})} ~\times ~  J_{n,\gothp, \cO} ~~
  \]
where $J_{n, \gothp, \cO}$ is the direct product of all of the $\Delta_{\gothp}$-components of $I_{n, \gothp, \cO} $ apart from the $\varepsilon_{\gothp}$-component. 
In terms of this decomposition,   $\xi$ is  in the image of 
$S_{\chi, \varepsilon}(\QQ_{n})$ under the restriction map if and only if $\xi_{\cO} \big(  J_{n,\gothp, \cO}\big) = 0$ for one (and hence for all) primes $\gothp$ lying over $p$.  Thus, the Selmer condition for $\xi$ is that $\xi_{\cO}$ factors through the quotient  $\Gal(M_n/K_n)_{\cO} \big/J_{n,p, \cO}$,   where $J_{n,p, \cO}$ denotes the $\cO$-submodule of    
 $\Gal(M_n/K_n)_{\cO}$ generated by all the  $J_{n,\gothp, \cO}$'s for $\gothp | p$.  
\smallskip

 For each $\gothp | p$, we have  $\cU_{n, \gothp, \cO} =  \cU_{n,\gothp, \cO}^{(\varepsilon_{\gothp})} ~ \times ~   \cV_{n,\gothp, \cO}$, where the second factor is the direct product of the $\Delta_{\gothp}$-components of $\cU_{n, \gothp, \cO} $ apart from the $\varepsilon_{\gothp}$-component. We then have
 $\cU_{n,p, \cO} ~=~  \cU_{n,p, \cO}^{[\varepsilon]} ~ \times ~   \cV_{n,p, \cO}$, where  we use the following notation:
\[
\cU_{n,p,\cO}^{[\varepsilon]} = \prod_{\gothp | p} ~ \cU_{n, \gothp, \cO}^{(\varepsilon_{\gothp})}~, \quad  \quad \quad \quad  
\cV_{n,p, \cO} = \prod_{\gothp | p} ~ \cV_{n, \gothp, \cO}~~
\]
The projection map from $\cU_{n,p, \cO}$ to  $\cU_{n,p, \cO}^{[\varepsilon]}$ will be denoted by $\pi_{n,\varepsilon}$.  We should point out that $\Delta_n$ acts naturally on $\cU_{n,p,\cO}$ and that both  $\cU_{n,p, \cO}^{[\varepsilon]}$ and  $\cV_{n,p, \cO}$ are invariant under the action of $\Delta_n$. The map  $\pi_{n,\varepsilon}$ is $\Delta_n$-equivariant and $\cO$-linear. 
\smallskip

Let $L_n$ denote the $p$-Hilbert class field of $K_n$,  the maximal, abelian, unramified $p$-extension of $K_n$.  We have $K_n \subseteq L_n \subset M_n$.  
Class field theory gives an exact sequence
\begin{equation}\label{eqn:cft}
1 \longrightarrow \overline{U}_{n,p} \longrightarrow   \cU_{n, p} \mathop{\longrightarrow}^{\alpha}   
 \Gal(M_n/K_n) \longrightarrow \Gal(L_n/K_n) \longrightarrow 1~~.
\end{equation}
The image of  $\cU_{n, \gothp}$ (as a factor in $\cU_{n,p}$) under the reciprocity map $\alpha$ is $I_{n,\gothp}$.  Tensoring the above exact sequence with $\cO$, the $\cO$-submodule $J_{n,p,\cO}$  of     $\Gal(M_n/K_n)_{\cO}$ is the image of 
$\cV_{n,p, \cO}$ under the map $\alpha_{_{\cO}}$.  The above observations then give us an exact sequence
\[
0 \longrightarrow  H^1_{unr}(\QQ_{n}, D) \longrightarrow S_{\chi, \varepsilon}(\QQ_{n}) \longrightarrow \Hom_{\cO[\Delta]}\Big( \cU_{n,p,\cO} \big/
 \cV_{n,p,\cO} \overline{U}_{n,p,\cO}, ~ D \Big) \longrightarrow 0~~.
\]
Note that  $\Gal(\QQ_n/\QQ)=\Gamma/\Gamma_n$ acts naturally on $H^1_{unr}(\QQ_{n}, D)$ and $S_{\chi, \varepsilon}(\QQ_{n})$.  Also, since $\Delta_n$ is identified with $\Delta \times \Gamma/\Gamma_n$,  there is also an action of $\Gamma/\Gamma_n$ on the next term of the above sequence. The maps are $\Gamma/\Gamma_n$-equivariant.

The restriction map  defines an isomorphism $H^1_{unr}(\QQ_n, D) \cong \Hom_{\Delta}(\Gal(L_n/K_n), ~D)$ and this group is finite.  We also have a canonical  isomorphism 
\[
\cU_{n,p,\cO}\Big/
 \cV_{n,p,\cO}\overline{U}_{n,p,\cO} ~ \cong ~ \cU_{n,p,\cO}^{[\varepsilon]} \Big/
 \pi_{n,\varepsilon}\big(\overline{U}_{n,p,\cO}\big)
\]
as   $\cO[\Delta_n]$-modules.   For brevity,  we will denote  $\pi_{n,\varepsilon}\big(\overline{U}_{n,p,\cO}\big)$ by  $\overline{U}_{n,p,\cO}^{[\ep]}$.  Thus,
\[
\overline{U}_{n,p,\cO}^{[\ep]} ~=~ \Im\big(\pi_{n,\varepsilon} \circ \lambda_{n,p,\cO}\big)~~.
\]
It is an $\cO[\Delta_n]$-submodule of $\cU_{n,p, \cO}^{[\varepsilon]}$.  
\medskip

In summary, $S_{\chi, \varepsilon}(\QQ_{n})$ contains the finite $\cO$-submodule  $H^1_{unr}(\QQ_{n}, D)$ and the corresponding quotient module is isomorphic to $\Hom_{\cO[\Delta]}\Big( \cU_{n,p, \cO}^{[\varepsilon]} \big/\overline{U}_{n,p,\cO}^{[\ep]}, ~ D\Big)$. 
In particular,  
$S_{\chi, \varepsilon}(\QQ_{n})$ is finite if and only if  $(\overline{U}_{n,p,\cO}^{[\ep]})^{(\chi)}$ has finite index in   $\big(\cU_{n,p, \cO}^{[\varepsilon]}\big)^{(\chi)}$.  We will verify this by showing that both $\cO$-modules have the same  $\cO$-rank.
\medskip

Let $\mu_{n,\gothp}$ denote the torsion subgroup of  $\cU_{n, \gothp}$.  The $p$-adic log map on  $\cU_{n, \gothp}$ is $\Gal(\cK_n/\QQ_p)$-equivariant and its image is open in the additive group of $\cK_n$.  Its kernel is     $\mu_{n,\gothp}$.  Thus, $\cU_{n, \gothp} \otimes_{\ZZ_p} \QQ_p$ is isomorphic to the regular representation of $\Gal(\cK_n/\QQ_p)$ over $\QQ_p$. Recall that $\Gal(\cK_n/\QQ_p)$ is identified with $\Delta_{\gothp} \times \Gamma/\Gamma_n$. Since $\varepsilon_{\gothp}$ is 1-dimensional,  it follows that   $\cU_{n, \gothp, \cO}^{(\varepsilon_{\gothp})} \otimes_{\cO} \cF$ is isomorphic to the regular representation of $\Gamma/\Gamma_n$ over $\cF$. In particular, one sees that the $\rank_{\cO}\big(\cU_{n, \gothp, \cO}^{(\varepsilon_{\gothp})}\big) = p^n$.  
$~$

A subscript $\cF$ will indicate tensoring over $\cO$ with $\cF$.  For example, we denote  $\cU_{n, \gothp, \cO} \otimes_{\cO} \cF$ by  $\cU_{n, \gothp, \cF}$. With this notation,  we have
\[
 \cU_{n, p, \cF}^{[\varepsilon]} ~=~  \prod_{\gothp | p}~  \cU_{n, \gothp, \cF}^{(\varepsilon_{\gothp})}~~.
\]
As a representation space for $\Delta_n=\Delta \times \Gamma/\Gamma_n$ over $\cF$, one can view $\cU_{n, p, \cF}^{[\varepsilon]}$ as an induced representation.  For this purpose, we single out one choice of $\gothp|p$ as in the introduction.  The decomposition subgroup of $\Delta_n$ for $\gothp$ is $\Delta_{\gothp} \times \Gamma/\Gamma_n$, which we denote by $\Delta_{n, \gothp}$. Furthermore, $\cU_{n, \gothp, \cF}^{(\varepsilon_{\gothp})}$ is a representation space for $\Delta_{n, \gothp}$ of dimension $p^n$ and is isomorphic to $\varepsilon_{\gothp} \otimes_{\cF} \sigma_{n}$, where $\sigma_n$ denotes the regular representation of $\Gamma/\Gamma_n$ over $\cF$.  We denote this representation of $\Delta_{n,\gothp}$ more briefly by $\varepsilon_{\gothp} \sigma_n$. We can also regard $\sigma_n$ as a representation of $\Delta_n$ over $\cF$ factoring through $\Delta_n/\Delta$. Then we have 
\[
\cU_{n, p, \cF}^{[\varepsilon]} ~\cong ~ \Ind_{\Delta_{n,\gothp}}^{\Delta_n}  \big( \varepsilon_{\gothp} \sigma_n \big) ~\cong ~  \Ind_{\Delta_{n,\gothp}}^{\Delta_n} \big(\varepsilon_{\gothp}\big) \otimes_{\cF} \sigma_n~~.
\]
One can regard $\Ind_{\Delta_{n,\gothp}}^{\Delta_n} \big(\varepsilon_{\gothp}\big)$ as the representation  $\Ind_{\Delta_{\gothp}}^{\Delta} \big(\varepsilon_{\gothp}\big)$ composed with the restriction map $\Delta_n \to \Delta$.

Since $\epp$ occurs with multiplicity 1 in $\chi |_{\Delta_{\gothp}}$, Frobenius Reciprocity implies that $\rho$ occurs with multiplicity 1 in  $\Ind_{\Delta_{\gothp}}^{\Delta} \big(\varepsilon_{\gothp}\big)$.   Now $\cF$ is generated over $\QQ_p$ by roots of unity of order prime to $p$ and hence is unramified.   The irreducible representations of $\Gamma/\Gamma_n$ over $\cF$ are $\tau_j$ for $0 \le j \le n$, where $\tau_j$ has kernel $\Gamma_j/\Gamma_n$. Then $\tau_j$ has degree $d_j=p^j - p^{j-1}$ for $1 \le j \le n$ and degree $d_0=1$ for $j=0$. We have 
\[
\sigma_n  ~\cong~ \bigoplus_{j=0}^{n} ~ \tau_j
\]
which has degree $p^n$. These remarks imply that each of the irreducible representations $\rho \otimes \tau_j$ occurs in $\cU_{n, p, \cF}^{[\varepsilon]} $ with multiplicity 1.  These are precisely the irreducible constituents in  $\big(\cU_{n, p, \cF}^{[\varepsilon]} \big)^{(\chi)}$ which therefore has $\cF$-dimension $dp^n$.  In contrast,  recall that  $\cU_{n, p, \cF}$ is isomorphic to the regular representation of $\Delta_n$ over $\cF$. The multiplicity of each $\rho \otimes \tau_j$ in $\cU_{n, p, \cF}$ is $d$ and therefore $\big(\cU_{n, p, \cF} \big)^{(\chi)}$ has $\cF$-dimension $d^2p^n$.

Now consider the $\Delta_n$-representation space $U_{n, \cF} = U_n \otimes_{\ZZ}\cF$. The absolutely irreducible constituents $\varphi$ of each $\sigma_n$ are 1-dimensional. Since $\QQ_n$ is totally real, each such $\varphi$ is an even character. Hence $d^{+}(\chi \varphi)=1$ and hence $\chi \varphi$ occurs with multiplicity 1 in  $U_n \otimes_{\ZZ} \overline{\QQ}_p$.  It follows that $\chi \tau_j$ occurs with multiplicity 1 in $U_{n,\cF}$ for each $j$, $0 \le j \le n$. Let $e_{\chi \tau_j}$ be the idempotent for $\chi \tau_j$ in $\cF[\Delta_n]$ (and actually in the center of $\cF[\Delta_n]$).   Then  $e_{\chi \tau_j}$ does not annihilate $U_{n,\cF}$. Furthermore, if one restricts $\chi \tau_j$ to $\Delta$, one obtains a multiple of $\chi$. 
Suppose that $\gothp$ is a prime of $K$ lying above $p$. It follows that $\epp$ occurs as a constituent in $\chi \tau_j |_{\Delta_{\gothp}}$. 
Let $e_{\epp} \in \cF[\Delta_{\gothp}]$ be the corresponding idempotent. We can regard $e_{\epp}$ as an element of $\cF[\Delta_n]$. 
For $0 \le j \le n$, let $\theta_j= e_{\epp} e_{\chi \tau_j}$. Thus $\theta_j \in \cF[\Delta_n]$. Also, the coefficients of $\theta_j$ are in $\overline{\QQ}$.   Furthermore, the above remarks show that $\theta_j$ does not annihilate $U_{n,\cF}$.  According to remark 2.6, it follows that the projection of $\Im(\lambda_{n,p,\cF})^{\theta_j}$ to $\cU_{n,\gothp.\cF}$ is nontrivial.  Note that the image of that projection is actually contained in $\cU_{n,\gothp,\cF}^{(\epp)}$.

 The above remarks show that, for $0 \le j \le n$, there exists an element $\alpha_{n,j}$ of $U_n$ such that $\beta_{n,j} = \lambda_{n,p}(\alpha_n)^{ \theta_j}$ has the following properties. Clearly, $\beta_{j,n}$ is in 
$\cU_{n,p,\cF}^{(\chi)}$ for the action of $\Delta$ and, more precisely, in the $\chi \tau_j$ component of $\cU_{n,p,\cF}$ for the action of $\Delta_n$.  Also, the image of $\beta_{j,n}$ under the composite map 
\[
\cU_{n,p,\cF} ~\longrightarrow~  \cU_{n,p,\cF}^{(\ep)} ~\longrightarrow ~\cU_{n,\gothp,\cF}^{(\epp)}
\]
is nontrivial. Here the first map is $\pi_{n, \ep, \cF}$ and the second map is the obvious projection map. 
It follows that $\pi_{n,\ep,\cF}(\beta_{j,n})$ is nontrivial. It must be in the $\chi \tau_j$-component of  $\cU_{n,p,\cF}^{(\ep)}$. Since $\chi \tau_j$ has multiplicity 1 in  $\cU_{n,p,\cF}^{(\ep)}$, the image of the $\chi \tau_j$-component of $U_{n,\cF}$ under the map $\pi_{n, \ep, \cF} \circ \lambda_{n,p,\cF}$ must coincide with the $\chi \tau_j$-component of 
 $\cU_{n,p,\cF}^{(\ep)}$.  Since this is so for all $j$, it follows that the image of $U_{n,\cF}^{(\chi)}$ under that map is precisely  $\big(\cU_{n,p,\cF}^{(\ep)}\big)^{(\chi)}$. This proves that 
 the image of $U_{n,\cO}^{(\chi)}$ under  $\pi_{n, \ep} \circ \lambda_{n,p,\cO}$ has the same $\cO$-rank as $\big(\cU_{n,p,\cO}^{(\ep)}\big)^{(\chi)}$ and hence has finite index. That means that $(\overline{U}_{n,p,\cO}^{[\ep]})^{(\chi)}$ indeed has finite index in   $\big(\cU_{n,p, \cO}^{[\varepsilon]}\big)^{(\chi)}$.      \hfill$\blacksquare$

 \bigskip
 
 One of the key observations in the proof of proposition 3.1 is the following result concerning the structure of $S_{\chi, \ep}(\QQ_{n})$ as an $\cO[\Gamma\big/\Gamma_n]$-module.
 \bigskip

\noindent {\bf Proposition 3.2.} {\em We have the following exact sequence of finite $\cO[\Gamma/\Gamma_n]$-modules:
\[
0 \longrightarrow  H^1_{unr}(\QQ_{n}, D) \longrightarrow S_{\chi, \varepsilon}(\QQ_{n}) \longrightarrow \Hom_{\cO[\Delta]}\Big( \big(\cU_{n,p, \cO}^{[\varepsilon]}\big)^{(\chi)} \big/ (\overline{U}_{n,p,\cO}^{[\ep]})^{(\chi)}, ~ D \Big) \longrightarrow 0~~.
\]}
\bigskip

\noindent 
In the classical case where $d=d^{+}=1$, $\chi$ is just an even character of $\Gal(K/\QQ)$ of order prime to $p$, where $K$ is a cyclic 
extension of $\QQ$. One can omit the $\ep$ everywhere in this case because $\chi$ determines $\ep$. Note that $\chi$ does not occur as a constituent in $\cV_{n,p,\cF}$ when $d=1$. The exact sequence in proposition 3.2 follows immediately from (\ref{eqn:cft}) by tensoring with $\cO$ and taking the $\chi$-component of each term. 
\bigskip

\section{ Selmer groups for Artin representations 2}  We will prove theorem 2 in this section. Here $\QQ_{\infty} = \bigcup_{n \ge 0} \QQ_n$ is the cyclotomic $\ZZ_p$-extension of $\QQ$ and $\QQ_n$ is the subfield of degree $p^n$ over $\QQ$. We can define the Selmer group either as a direct limit or as the kernel of a global-to-local map. That is, 
\[
S_{\chi, \varepsilon}(\QQ_{\infty}) ~=~ \varinjlim_n ~S_{\chi, \ep,}(\QQ_n) ~=~ \ker\Big(H^1(\QSQ_{\infty}, \cD) \longrightarrow  \prod_{\ell \in \Sigma} ~ \cH_{\ell}(\QQ_{\infty}, D)\Big)
\]
where the direct limit is induced by the restriction maps $H^1(\QSQ_n, D) \to H^1(\QSQ_m, D)$ for $m \ge n \ge 0$.  Note that, for $\ell \neq p$, $\Gal(\QQ_{\infty, \nu}^{unr}\big/\QQ_{\infty, \nu})$ has profinite order prime to $p$ and so an element in $H^1(\QSQ_{\infty}, D)$ is locally unramified at $ \nu | \ell$ if and only if it is locally trivial. Thus, we can define 
\[
\cH_{\ell}(\QQ_{\infty}, D) ~= ~\prod_{\nu | \ell} ~ H^1(\QQ_{{\infty}, \nu}, D) ~~
\]
for $\ell \neq p$.  For $\ell=p$, one defines  
\[
 \cH_{p}(\QQ_{\infty}, D)  = \Im \left(   H^1(\QQ_{\infty, \pi}, D) \longrightarrow H^1(\QQ_{\infty, \pi}^{unr},   D/D^{(\varepsilon_{\gothp})})\right)
 = H^1(\QQ_{\infty, \pi}, D)/L(\QQ_{\infty, \pi}, D)~ 
\]
where $L(\QQ_{\infty, \pi}, D)$ is a $\LO$-submodule of   $H^1(\QQ_{\infty, \pi},D)$ which sits in the exact sequence
\[
0~ \longrightarrow ~  H^1(\QQ_{\infty, \pi}, D^{(\epp)} ) ~ \longrightarrow ~ L(\QQ_{\infty, \pi}, D) ~ \longrightarrow ~ H^1_{unr}(\QQ_{\infty,\pi}, D') \longrightarrow ~0~~.
\]
Here we write $D'$ for $D/D^{(\epp)}$ for brevity and  $H^1_{unr}(\QQ_{\infty,\pi}, D')$ denotes the kernel of the restriction map 
$H^1(\QQ_{\infty,\pi}, D') \to  H^1(\QQ^{unr}_{\infty,\pi}, D')$.  Also, since $p \nmid |\DP|$, we have $D \cong D^{(\epp)} \oplus D'$. This implies the injectivity in the above sequence and also implies that the above sequence splits. Thus,  $L(\QQ_{\infty, \pi}, D)$ can be identified with
 \[
 H^1(\QQ_{\infty, \pi}, D^{(\epp)} ) \oplus H^1_{unr}(\QQ_{\infty,\pi}, D')~~,  
 \]
 considered in the obvious way as a $\LO$-submodule of  $H^1(\QQ_{\infty, \pi}, D)$.
\smallskip

 Our study of   $S_{\chi, \varepsilon}(\QQ_{\infty})$ will be based on the following control theorem.  
\bigskip

\noindent {\bf Proposition 4.1.} {\em Assume that Hypothesis A is satisfied and that $\chi$ is not the trivial character.   For any $n \ge 0$, we have an exact sequence of $\cO[\Gamma/\Gamma_n]$-modules:
\[
0 ~\longrightarrow ~ S_{\chi, \varepsilon}(\QQ_{n})~ \longrightarrow ~ S_{\chi, \varepsilon}(\QQ_{\infty})^{\Gamma_n} ~ \longrightarrow ~ H^0(\Delta_{\gothp}, ~D/D^{(\epp)}) ~ \longrightarrow ~ 0~~.
\]
Note that  $H^0(\Delta_{\gothp}, D/D^{(\epp)}) \cong \big(\cF/\cO\big)^t$, where $t$ is the multiplicity of the trivial representation of $\Delta_{\gothp}$  in $V/V^{(\epp)}$ and the action of $\Gamma/\Gamma_n$ of $\big(\cF/\cO\big)^t$ is trivial. }
\bigskip

\noindent {\em Proof.}   We exclude the case where $\chi$ is the trivial character. We apply the snake lemma to the following commutative diagram.
\[
\xymatrix {
  0  \ar[r] & S_{\chi, \ep}(\QQ_n)   \ar[r]  \ar[d]^{s_n} & H^1(\QQ_{\Sigma}/\QQ_n, ~D) \ar[r]^{\phi_n}  \ar[d]^{h_n} &  \prod_{\ell \in \Sigma} ~ \cH_{\ell}(\QQ_{n}, D)  \ar[d]^{r_n} \\
0 \ar[r] & S_{\chi, \ep}(\QQ_{\infty}, D)^{\Gamma_n}  \ar[r]  & H^1(\QQ_{\Sigma}/\QQ_{\infty}, \cD)^{\Gamma_n}  \ar[r] &   \prod_{\ell \in \Sigma} ~ \cH_{\ell}(\QQ_{\infty}, D)^{\Gamma_n}  }
\]
 where the vertical maps are the obvious restriction maps and $\phi_n$ is the global-to-local map defining $S_{\chi, \ep}(\QQ_n)$.   
 
 We first show that the $\cO$-corank of $H^1(\QQ_{\Sigma}/\QQ_n, ~D)$ is equal to $d^{-}p^n=(d-1)p^n$. To see this, we have an isomorphism  $H^1(\QQ_{\Sigma}/\QQ_n, ~D) \cong  H^1 (\QQ_{\Sigma}/\QQ, ~\Ind_{\QQ_n}^{\QQ}(D))$ by Shapiro's Lemma. . One can identify 
 $\Ind_{\QQ_n}^{\QQ}(D)$ with $D \otimes_{\cO}\cO[\Gamma/\Gamma_n]$ and its $\cO$-corank is $dp^n$.  
  Let $\cF_n=\cF(\mu_{p^n})$ and let $\cO_n$ denote the ring of integers in that local field.  Tensoring with $\cO_n$ over $\cO$, it suffices to show that $H^1\big(\QSQ, D\otimes_{\cO}\cO_n[\Gamma/\Gamma_n]\big)$ has $\cO_n$-corank equal to $d^-p^n$.   
 Now $\Ind_{\QQ_n}^{\QQ}(V) \otimes_{\cF} \cF_n$ is isomorphic to the direct sum of the irreducible representations $\rho \otimes \varphi$ over $\cF_n$, where $\varphi$ varies over the characters of $\Gamma/\Gamma_n$, all of which are $d$-dimensional and have $d^{+}=1, ~ d^{-}=d-1$.  For brevity, we let $V_{\varphi}$ and $T_{\varphi}$ denote the corresponding $\cF_n$-representation spaces and Galois invariant $\cO_n$-lattices.  Let $D_{\varphi}=V_{\varphi}/T_{\varphi}$.   Propositions 2.1 and 2.2 and corollary 2.7 imply that the $\cO_n$-corank of $H^1(\QSQ, D_{\varphi})$ is $d^{-}$ for all $\varphi$.  Hence the $\cO_n$-corank of  $H^1(\QQ_{\Sigma}/\QQ, \oplus_{\varphi} D_{\varphi})$ is $d^{-}p^n$. There is a $\Delta$-equivariant surjective homomorphism from $D \otimes_{\cO} \cO_n[\Gamma/\Gamma_n]$ to $\oplus_{\varphi} D_{\varphi}$ with finite kernel. It follows that the $\cO_n$-corank in question is indeed equal to $d^-p^n$. 
 
 We now show that $\cH_p(\QQ_n,D)$ also has $\cO$-corank equal to $d^{-}p^n$ and that $\cH_{\ell}(\QQ_n, D)$ is finite for all $\ell \in \Sigma, ~\ell \neq p$. We use the Poitou-Tate formulas for the local Euler-Poincar\'e characteristics . Let $T^{*}$ denote  the compact Galois module $\Hom(D, \mu_{p^{\infty}})$.   It is free as an $\cO$-module. 
 For any prime $\ell$, let $\cL$ be a finite extension of $\QQ_{\ell}$. Then Poitou-Tate duality implies that $H^2(\cL, D)$ is dual to $H^0(\cL, T^{*})$.  Since the action of $G_{\cL}$ on $D$ factors through a finite quotient and through an infinite quotient on $\mu_{p^{\infty}}$, it follows that  $H^0(\cL, T^{*})=0$ and hence $H^2(\cL, D)=0$. The only assumption we need is that $G_{\cL}$ acts on $D$ through a finite quotient group. 
 
As before, let $D'$ denote $D/D^{(\epp)}$. Then $\corank_{\cO}(D')=d-1=d^{-}$.   The Euler-Poincar\'e characteristic formulas for $D'$ imply that
 \[
 \corank_{\cO}\big( H^1(\QQ_{n,\pi},  D' )\big) = d^{-}p^n + \corank_{\cO}\big( H^0(\QQ_{n,\pi}, D') \big)
 \]
 where $\pi$ is the unique prime above $p$ in $\QQ_n$ (or $\QQ_{\infty}$ as before). Here we are using the vanishing of $H^2(\QQ_{n,\pi}, D')$.  Since $\G_{\QQ_p}$ acts on $D'$ through a finite quotient of order prime to $p$, it follows that   
\[
H^0(\QQ_{n,\pi}, D') ~=~H^0(\QQ_{p}, D') ~\cong ~(\cF/\cO)^t
\]
which is $\cO$-divisible and has $\cO$-corank $t$.  However, letting $\QQ_{n,\pi}^{unr,p}$ denote the unramified $\ZZ_p$-extension of $\QQ_{n,\pi}$, the inflation-restriction sequence shows that 
 \[
 \ker\big( H^1(\QQ_{n,\pi}, D') \to H^1(\QQ_{n,\pi}^{unr}, D')\big) ~=~  \ker\big( H^1(\QQ_{n,\pi}, D') \to H^1(\QQ_{n,\pi}^{unr,p}, D')\big) 
 \]
 \[
 \cong~ H^1\big(\QQ_{n,\pi}^{unr,p}/\QQ_{n,\pi}, ~H^0(\QQ_{n,\pi}^{unr,p}, D') \big)~=~ H^1\big(\QQ_{n,\pi}^{unr,p}/\QQ_{n,\pi},~ H^0(\QQ_p,  D')\big)~~, 
 \]
 where the last equality follows again because the Galois action on $D'$ is through a quotient group of order prime to $p$. It follows that this kernel has $\cO$-corank $t$ and therefore that $\cH_p(\QQ_n,~D)$ indeed has $\cO$-corank equal to $d^{-}p^n$. 
 
 Now suppose $\ell \neq p$ and let $\nu$ be a prime of $\QQ_n$ above $\ell$.  In this case, we have the formula
 \[
 \corank_{\cO}\big( H^1(\QQ_{n,\nu}, D )\big) =  \corank_{\cO}\big( H^0(\QQ_{n,\nu}, D)\big) .
 \]
This $\cO$-corank could be positive. But, just as above, the inflation-restriction sequence shows that the $\cO$-corank of the kernel of the restriction map $H^1(\QQ_{n,\nu}, D ) \to H^1(\QQ^{unr}_{n,\nu}, D )$ is the same as the $\cO$-corank of  $H^0(\QQ_{n,\nu}, D)$ and hence the image of the restriction map is finite.  The finiteness of $\cH_{\ell}(\QQ_n, D)$ 
follows.

 Thus, both $H^1(\QSQ_n,~D)$ and $\prod_{\ell \in \Sigma}~ \cH_{\ell}(\QQ_n,  D)$ have the same $\cO$-corank.  Furthermore, the kernel of the global-to-local map $\phi_n$ is $S_{\chi,\ep}(\QQ_n)$ which is finite by proposition 3.1.  Therefore, the cokernel of $\phi_n$ has $\cO$-corank 0 and must be finite.  This is one of the hypotheses in proposition 3.2.1 in \cite{gre10}. Another hypothesis is called {\bf LEO}$(D)$ there, but that follows immediately from corollary 2.7 and proposition 2.2. A third hypothesis is that no subquotient of $D[p] \cong T/pT$ is isomorphic to $\mu_p$ as a $G_{\QQ_n}$-module. However, this is satisfied because $\Delta$ has order prime to $p$, hence the reduction of $\rho$ modulo $p$ is still irreducible as a representation space over the residue field $\cO/(p)$, and has no subquotient which is odd. Consequently, proposition 3.2.1 in \cite{gre10} implies that $\phi_n$ is surjective. We can now apply the snake lemma to the commutative diagram at the beginning of this proof.
 
 First of all, $\ker(h_n) \cong H^1\big(\Gamma_n, H^0(\QQ_{\infty}, D)\big)$.  Since the action of $G_{\QQ_n}$ on $D$ factors through a quotient of order prime to $p$ and $\chi$ is nontrivial and irreducible, it follows that  $H^0(\QQ_{\infty}, D)= H^0(\QQ, D)=0$.  Therefore, $h_n$ is injective, and so is $s_n$.  Also, since $\Gamma_n \cong \ZZ_p$, it has $p$-cohomological dimension 1. This implies that $\coker(h_n)=0$.  The snake lemma implies that $\coker(s_n) \cong \ker(r_n)$.  For each $\ell \in \Sigma$, consider the restriction map
 \[
 r_{n, \ell}: ~ \cH_{\ell}(\QQ_n, D) ~ \longrightarrow ~ \cH_{\ell}(\QQ_{\infty}, D)~~.
 \]
  Now, if $\ell \neq p$, then $\ell$ is unramified in $\QQ_{\infty}/\QQ$. It follows that, for $\nu | \ell$, $\QQ_{n, \nu}^{unr} = \QQ_{\infty,\nu}^{unr}$. This implies that $r_{n, \ell}$ is injective.  To determine $\ker(r_{n,p})$, let us write $\Gamma_n^{unr}$ for $\Gal(\QQ_{n,\pi}^{unr,p}/\QQ_{n,\pi})$ and $\Gamma_n^{cyc}$ for $\Gamma_n=\Gal(\QQ_{\infty,\pi}/\QQ_{n,\pi})$. Just as before, we have 
 \[
 \ker\big(H^1(\QQ_{n,p}, D') \to H^1(\QQ_{\infty,\pi} \QQ_{n,\pi}^{unr}, D')\big) ~=~ \ker\big(H^1(\QQ_{n,p}, D') \to H^1(\QQ_{\infty,\pi} \QQ_{n,\pi}^{unr,p}, D') \big)
  \]
  \[ 
  \cong ~ H^1\big(\QQ_{\infty,\pi} \QQ_{n,\pi}^{unr,p}\big/\QQ_{n, \pi}, ~H^0(\QQ_p, D')\big) ~\cong ~ \Hom\big(\Gamma_n^{cyc} \times \Gamma_n^{unr}, ~H^0(\QQ_p, D')\big)~~
  \]
  \[
  \cong ~  \Hom\big(\Gamma_n^{cyc}, ~H^0(\QQ_p, D')\big) \times  \Hom\big(\Gamma_n^{unr}, ~H^0(\QQ_p, D')\big)~.
  \]
  The second factor is the kernel of the restriction map $H^1(\QQ_{n,\pi}, D') \to H^1(\QQ_{n,\pi}^{unr}, D')$ and therefore $\ker(r_{n, p})$ 
  is isomorphic to the first factor, and hence isomorphic (non-canonically) to $H^0(\QQ_p, D') = H^0(\DP, ~D')$. 
    It is a cofree $\cO$-module of corank $t$ and $\Gamma/\Gamma_n$ acts trivially. We have proved that $s_n$ is indeed injective and has cokernel as stated in the proposition.   \hfill$\blacksquare$

\bigskip
  
 As always, we exclude in the above proposition the case where $\chi$ is the trivial character. In that case, one easily shows that $S_{\chi, \ep}(\QQ_{\infty})=0$. This is so even though $S_{\chi, \ep}(\QQ_n)$ is infinite for all $n \ge 0$ according to (\ref{eqn:trivial}).  And so, $\ker(s_n)$ is infinite when $\chi$ is the trivial character.

\bigskip

The next result is theorem 2 in the introduction. 
 \bigskip
 
\noindent {\bf Proposition 4.2.}  {\em Suppose that  Hypothesis A is satisfied.    Then  $S_{\chi, \varepsilon}(\QQ_{\infty})$ is a cofinitely generated, cotorsion $\Lambda_{\cO}$-module. }
\bigskip

\noindent {\em Proof.} By definition,  $S_{\chi, \varepsilon}(\QQ_{\infty})$ is a $\Lambda_{\cO}$-submodule of $H^1(\QQ_{\Sigma}/\QQ_{\infty}, D)$ and hence its Pontryagin dual $X_{\chi, \varepsilon}$ is a quotient of the Pontryagin dual of $H^1(\QQ_{\Sigma}/\QQ_{\infty}, D)$ as a $\Lambda_{\cO}$-module.  
Proposition 3 in \cite{gre89} asserts that $H^1(\QQ_{\Sigma}/\QQ_{\infty}, D)$ is cofinitely generated as a $\Lambda$-module, where $\Lambda=\ZZ_p[[\Gamma]]$.  Hence its Pontryagin dual is finitely generated as a $\Lambda$-module and hence as a $\Lambda_{\cO}$-module. Therefore, $X_{\chi, \varepsilon}$ is also 
a finitely generated $\Lambda_{\cO}$-module.  

 If $X$ is a finitely generated $\Lambda_{\cO}$-module of rank $r$, then 
 \[
 \rank_{\cO} (X_{\Gamma_n}) =rp^n + O(1) 
 \] 
 as $n \to \infty$. If $S$ is the Pontryagin dual of $X$, then $\corank_{\cO}(S^{\Gamma_n})=\rank(X_{\Gamma_n})$.  Thus, $r=0$ if and only if $\corank_{\cO}(S^{\Gamma_n})$ is bounded. Considering $S=S_{\chi,\ep}(\QQ_{\infty})$,  propositions 3.1 and 4.1 imply the boundedness of these $\cO$-coranks. Hence $r=0$ and that means that $S_{\chi, \ep}(\QQ_{\infty})$ is $\Lambda_{\cO}$-cotorsion.    \hfill$\blacksquare$
\bigskip

\noindent {\bf Remark 4.3.} This remark concerns so-called {\em trivial zeros} or {\em exceptional zeros}.  Suppose that the trivial character of $\Delta_{\gothp}$ occurs as a constituent in $\chi |_{\Delta_{\gothp}}$ with multiplicity $t$, but that $\epp$ is nontrivial.  Proposition 4.1 then implies that $S_{\chi, \varepsilon}(\QQ_{\infty})^{\Gamma}$ is infinite and that its $\cO$-corank is equal to $t$.  Thus, $(X_{\chi, \ep})_{\Gamma}$ has $\cO$-rank $t$.  As in the introduction, let $\theta_{\chi, \ep}$ be a generator of the characteristic ideal of the $\Lambda_{\cO}$-module $X_{\chi, \ep}$. If $\varphi: \Gamma \to \overline{\QQ}^{\times}_p$ is a continuous group homomorphism, then $\varphi$ can be extended uniquely to a continuous $\cO$-algebra homomorphism $\varphi: \Lambda_{\cO} \to \overline{\QQ}_p$. In particular, if $\varphi_0$ is the trivial character of $\Gamma$, then one obtains a continuous surjective $\cO$-algebra  homomorphism $\varphi_0: \Lambda_{\cO} \to \cO$ whose kernel is generated by $\gamma_0-1$. Here $\gamma_0$ is a fixed topological generator of $\Gamma$.  The fact that $(X_{\chi, \ep})_{\Gamma}$ has $\cO$-rank $t$ implies that $(\gamma_0-1)^t$ divides $\theta_{\chi, \ep}$  in $\Lambda_{\cO}$.  Thus, $\varphi_0(\theta_{\chi,\ep})=0$. One would say that $\varphi_0$ is a zero of $\theta_{\chi, \ep}$ of order $\ge t$. It seems reasonable to conjecture that the order of vanishing is exactly $t$.

     As usual in Iwasawa theory, the $\cO$-algebra $\Lambda_{\cO}$ is isomorphic to the formal power series ring $\cO[[T]]$. One defines such an isomorphism by sending 
 $\gamma_0-1$  to $T$. Thus, $\theta_{\chi, \ep}$ corresponds to a power series $f_{\chi, \ep}(T)$ in $\cO[[T]]$. Under the above assumption about $\chi |_{\Delta_{\epp}}$, it follows that 
 $f_{\chi, \ep}(T)=T^t g_{\chi, \ep}(T)$, where $g_{\chi, \ep}(T) \in \cO[[T]]$. The conjecture is that $g_{\chi, \ep}(0) \neq 0$. 
 
 One interesting example is discussed in detail in \cite{gre91}, pages 227-231. The representation $\rho$ in that example is the 2-dimensional irreducible representation of $\Delta=\Gal(K/\QQ)$, where $K$ is the splitting field over $\QQ$ for $x^3-x+1$. In fact, $\Delta \cong S_3$ and $\rho$ is realizable over $\QQ$. We can take $\cF=\QQ_p$ and $\cO=\ZZ_p$.  We take $p=23$. Then $\Delta_{\gothp}$ is a subgroup of order 2. In \cite{gre91}, the case where $\epp$ is the nontrivial character of $\Delta_{\gothp}$ is discussed. Then $\Delta_{\gothp}$ acts trivially on $V/V^{(\epp)}$ and $t=1$. It is shown there that $S_{\chi, \ep}(\QQ_{\infty}) \cong \QQ_p/\ZZ_p$. Of course, $\Gamma$ acts trivially. And so, the characteristic ideal of $X_{\chi, \ep}$ is generated by $T$.  The calculation in \cite{gre91} is rather subtle and depends on calculating something which one could call the $\cL$-invariant for the representation $\rho$ and $p=23$.  In particular, the $\cL$-invariant turns out to be nonzero. It would be tempting to extend such a calculation to more general $\rho$'s satisfying Hypothesis A, but we have done nothing in that direction. 
\bigskip

\noindent {\bf Remark 4.4.} One can ask whether $H_{unr}^1(\QQ_{\infty}, D)^{\Gamma}$ can be infinite. Let $t$ be as in remark 4.3. If $t=0$, then $S_{\chi, \varepsilon}(\QQ_{\infty})^{\Gamma}$ is finite and hence so is $H_{unr}^1(\QQ_{\infty}, D)^{\Gamma}$. Suppose now that $t=1$.  Let $\ep_0$ denote the trivial character of $\Gal(\cK/\QQ_p)$. Hypothesis A is then satisfied and hence $S_{\chi, \ep_0}(\QQ_{\infty})$ is $\LO$-cotorsion by proposition 4.2.  
Actually, proposition 4.1  implies that $S_{\chi, \ep_0}(\QQ_{\infty})^{\Gamma}$ is finite. Since $H_{unr}^1(\QQ_{\infty}, D)  \subseteq S_{\chi, \ep_0}(\QQ_{\infty})$,  it follows that $H_{unr}^1(\QQ_{\infty}, D)^{\Gamma}$ is finite in that case too. However, if $t \ge 2$, then it turns out that $H_{unr}^1(\QQ_{\infty}, D)^{\Gamma}$ is infinite. Its $\cO$-corank will be at least $t-1$. We will discuss this in the sequel to this paper.

 \bigskip

 Another basic result concerning  $S_{\chi, \varepsilon}(\QQ_{\infty})$ is the following {\em purity} result concerning the $\Lambda_{\cO}$-module $X_{\chi, \ep}$ and asserts equivalently that  $S_{\chi, \varepsilon}(\QQ_{\infty})$ is an {\em almost divisible} $\Lambda_{\cO}$-module in the sense of \cite{gre06}. 
\bigskip

\noindent {\bf Proposition 4.5.} {\em Assume that Hypothesis A is satisfied. Then the Pontryagin dual $X_{\chi, \ep}(\QQ_{\infty})$  of $S_{\chi,\varepsilon}(\QQ_{\infty})$  has no nontrivial, finite $\Lambda_{\cO}$-submodules. }
\bigskip

\noindent {\em Proof.}   This result follows from 4.1.1 in \cite{gre16}.  However, that proposition  is formulated in terms of a Galois representation over $\Lambda_{\cO}$ instead of over $\cO$, namely the so-called cyclotomic deformation of $\rho$ discussed in section 3 of \cite{gre94}. Proposition 3.2 in that paper provides the $\Lambda_{\cO}$-module isomorphism  between $S_{\chi, \ep}(\QQ_{\infty})$ and the Selmer group over $\QQ$ associated to the cyclotomic deformation of $\rho$.  It is not difficult to verify the hypotheses (some of which were already discussed in proving the surjectivity of $\phi_n$ in the proof of proposition 4.1.)  One hypothesis is  that  $L(\QQ_{\infty,\pi}, D)$ (which was defined above) is almost divisible as a $\Lambda_{\cO}$-module.    
This just means that $L(\QQ_{\infty,\pi}, D)$  has no proper $\LO$-submodules of finite index. Now $H^1_{unr}(\QQ_{\infty,\pi}, D')$ is isomorphic to $(\cF/\cO)^t$ and is divisible as an $\cO$-module. Furthermore, 
$H^1(\QQ_{\infty, \pi}, D^{(\epp)} )$ is also divisible as an $\cO$-module because $D^{(\epp)}$ is divisible and $G_{\QQ_{\infty,\pi}}$ has $p$-cohomological dimension 1. 
\hfill$\blacksquare$

\bigskip

\noindent {\bf Remark 4.6.}  The conclusion in proposition 4.5 means that $S_{\chi,\ep}(\QQ_{\infty})$ is an almost divisible $\LO$-module.  One variant of the Selmer group $S_{\chi, \ep}(\QQ_{\infty})$ is the so-called {\em strict} Selmer group $S^{str}_{\chi, \ep}(\QQ_{\infty})$ where we continue to require elements of $H^1(\QSQ_{\infty}, D)$ to be locally unramified at all $\nu | \ell$ for $\ell \neq p$, but for $\ell=p$, we require those elements to have trivial image in $H^1(\QQ_{\infty,\pi}, D')$ (rather than in $H^1(\QQ^{unr}_{\infty,\pi}, D')$).  However, this just means that we replace 
\[
L(\QQ_{\infty,\pi}, D)~=~H^1(\QQ_{\infty, \pi}, D^{(\epp)} ) \oplus H^1_{unr}(\QQ_{\infty,\pi}, D')
\]
by $L^{str}(\QQ_{\infty,\pi}, D)~=~H^1(\QQ_{\infty, \pi}, D^{(\epp)} )$. Since this is an almost divisible $\LO$-module, it again follows that it $S^{str}_{\chi, \ep}(\QQ_{\infty})$ is almost divisible as a $\LO$-module. In section 5, it will be useful to know the same thing for the variant where we replace $L(\QQ_{\infty,\pi}, D)$ by
\[
L'(\QQ_{\infty,\pi}, D)~=~ H^1(\QQ_{\infty, \pi}, D^{(\epp)} )_{_{\LO-div}} ~\oplus ~H^1_{unr}(\QQ_{\infty,\pi}, D')~~.
\]
Here, for a discrete $\LO$-module $H$, we let $H_{_{\LO-div}}$ denote the maximal $\LO$-divisible submodule of $H$.  As we point out in section 5, $H^1(\QQ_{\infty, \pi}, D^{(\epp)} )$ fails to be $\LO$-divisible just in the case where  $\ep =\omega$. We denote the corresponding Selmer group by  $S'_{\chi,\ep}(\QQ_{\infty})$. Since  the $\LO$-module $L'(\QQ_{\infty,\pi}, D)$ is almost divisible, it will again follow that $S'_{\chi,\ep}(\QQ_{\infty})$ is also an almost divisible $\LO$-module. 
\bigskip

 There is one case where one can prove the analogue of proposition 4.1 even if $d^{+} \neq 1$, namely the case where $d^{+}=d-1$. Then $d^{-}=1$. 
We use the same notation as before, i.e., $\rho, V, T, D,  \cF,  \cO$, etc..
Just as in Hypothesis A, we will assume that $V$ has a $\Delta_{\gothp}$-invariant subspace $W_{\gothp}$ such that $\dim_{\cF}(W) =d^{+}$ and that the $\Delta_{\gothp}$-representation spaces $W$ and $V/W$ have no irreducible constituents in common.  Now $\Delta_{\gothp}$ acts on $V/W$ by a character $\epp'$ and the assumption means that $\epp'$ has multiplicity 1 in $V$.  Let $\epp = \chi - \epp'$, which is the character of $W$ as a $\Delta_{\gothp}$-representation space. We will denote $W$ by $V^{(\epp)}$ even though $\epp$ is not necessarily an irreducible character of $\Delta_{\gothp}$. We use a similar notation for $D$. Thus, $D \cong D^{(\epp)} \oplus D^{(\epp')}$, where the two summands have $\cO$-coranks $d-1$ and $1$, respectively. With this notation, we can define $S_{\chi, \ep}(\QQ_{\infty})$ exactly as in the introduction.  It turns out that we can use proposition 4.2 (for a different Artin representation) to prove the following result.

\bigskip

\noindent {\bf Proposition 4.7.}  {\em Under the above assumptions,   $S_{\chi, \varepsilon}(\QQ_{\infty})$ is a cofinitely generated, cotorsion $\Lambda_{\cO}$-module. }
\bigskip

\noindent {\em Proof.}  The fact that the Selmer group is cofinitely generated as a $\Lambda_{\cO}$-module follows from the fact that the same is true for $H^1(\QQ_{\Sigma}/\QQ, ~D)$.
We can assume that $\mu_p \subset K$.  Let $\cF(\omega)$ be the 1-dimensional $\cF$-representation space where $\Delta$ acts by character $\omega$. Here $\omega$ gives the action of $\Delta$ on $\mu_p$ and its values are the $(p-1)$-st roots of unity (which are in $\cF$).  We let $U=\Hom\big(V, ~\cF(\omega)\big)$ which is  the underlying $\cF$-representation space for an Artin representation $\sigma$ factoring through $\Delta$. Let $\psi$ be the character of $\sigma$. Note that $d^{+}(\psi)=d^{-}(\chi)=1$. Furthermore, the orthogonal complement of $W$ is the $\Delta_{\gothp}$-invariant subspace $U^{(\dpp)}$ of $U$ where $\dpp = \omega \epp'^{-1}$.  We let $\delta$ denote the corresponding character of $\Gal(\cK/\QQ_p)$.  Both $\dim_{\cF}(U^{(\dpp)})$ and $d^{+}(\psi)$ are equal to 1.  Furthermore, $\delta_{\gothp}$ is not a constituent in $U/U^{(\dpp)}$.  And so we are in the situation considered previously. Proposition 3.4 implies that $S_{\psi, \delta}(\QQ_{\infty})$ is a cotorsion $\Lambda_{\cO}$-module.  We will use arguments in \cite{gre89} to deduce that  $S_{\chi, \varepsilon}(\QQ_{\infty})$ is also $\LO$-cotorsion. In fact, there is a close relationship between the structures of $S_{\psi, \delta}(\QQ_{\infty})$ and $S_{\chi, \ep}(\QQ_{\infty})$ as $\Lambda_{\cO}$-modules. 

The action of $G_{\QQ}$ on $\mu_{p^{\infty}}$ is given by $\omega \kappa$, where $\omega$ and $\kappa$ are as defined in the introduction. The character $\kappa$ factors through $\Gamma = \Gal(\QQ_{\infty}/\QQ)$. We have $D=V/T$ as usual.  We define $T^{*} = \Hom(D, \mu_{p^{\infty}})$, which is a free $\cO$-module of rank $d=\dim_{\cF}(V)$.  Let   $V^{*} = T^{*}\otimes_{\cO}\cF$, and let $D^{*} = V^{*}/T^{*}$.  The representation space  $V^{*}$ of $G_{\QQ}$ is not an Artin representation. However, if we let $\cF(\omega \kappa)$ and $\cF(\kappa)$ denote the 1-dimensional vector spaces over $\cF$ on which $G_{\QQ}$ acts by the cyclotomic character $\omega \kappa$ and by $\kappa$, respectively,  then we have
\[
V^{*} ~ \cong ~ \Hom_{\cF}\big(V, ~ \cF(\omega \kappa)\big)~ \cong~U \otimes_{\cF} \cF(\kappa) ~~.
\]
We will just write $V^{*} \cong U \otimes \kappa$ which we think of as the twist of the Artin representation $U$ by $\kappa$. We have a filtration  $0 \subset W \subset V$ and the orthogonal complements give a filtration on $U$ and on $U\otimes \kappa$.   Proposition 11 in \cite{gre89} deals with an analogous situation, although the filtrations used there are defined in a different way. The Selmer groups denoted by $S_{V/T}(\QQ_{\infty})$ and $S_{V^{*}/T^{*}}(\QQ_{\infty})$ in that paper are also defined in terms of the filtrations, completely analogously  to the definitions here. Proposition 11 implies that if one of those Selmer groups is $\Lambda_{\cO}$-cotorsion, then so is the other. However, $\kappa$ factors through $\Gamma=\Gal(\QQ_{\infty}/\QQ)$ and so $V^{*} \cong U$ as a representation space for $G_{\QQ_{\infty}}$. It would then follow that $S_{D^{*}}(\QQ_{\infty})$ and $S_{\psi, \delta}(\QQ_{\infty})$ are the same as sets. More precisely, for the action of $\Gamma$, one has $S_{D^{*}}(\QQ_{\infty}) \cong S_{\psi, \delta}(\QQ_{\infty})\otimes \kappa$. Therefore, the fact that  $S_{\psi, \delta}(\QQ_{\infty})$ is $\Lambda_{\cO}$-cotorsion implies the same for  $S_{D^{*}}(\QQ_{\infty})$ and that in turn implies that $S_D(\QQ_{\infty})=S_{\chi, \varepsilon}(\QQ_{\infty})$ is $\Lambda_{\cO}$-cotorsion.   \hfill$\blacksquare$

\bigskip

\noindent {\bf Remark 4.8.}  With the notation and assumptions as in proposition 4.7 and its proof, let $X_{\chi, \ep}(\QQ_{\infty}), ~ X_{D^*}(\QQ_{\infty})$, and $X_{\psi, \delta}(\QQ_{\infty})$ denote the Pontryagin duals of $S_{\chi, \ep}(\QQ_{\infty})=S_D(\QQ_{\infty})$, $~$ $S_{D^*}(\QQ_{\infty})$, and $S_{\psi, \delta}(\QQ_{\infty})$, respectively. They are finitely generated, torsion $\Lambda_{\cO}$-modules. Let $\theta_{\chi, \ep}, ~\theta_{D^*}$, and $\theta_{\psi, \delta}$ denote generators for the corresponding characteristic ideals in $\Lambda_{\cO}$.  The first two are related (up to a unit) by the involution $\iota$ of the $\cO$-algebra $\Lambda_{\cO}$ induced by the automorphism $\gamma \mapsto \gamma^{-1}$ of $\Gamma$. This is the content of  theorem 2 in \cite{gre89} whose proof applies without change to the case at hand. Thus, we could simply chose    $\theta_{\chi, \ep}$ to be $\theta_{D^*}^{\iota}$. On the other hand, we have $X_{D^*}(\QQ_{\infty}) \cong X_{\psi, \delta}(\QQ_{\infty}) \otimes \kappa^{-1}$.  Thus, 
\[
X_{\psi, \delta}(\QQ_{\infty}) ~ \cong ~ X_{D^*}(\QQ_{\infty})\otimes \kappa 
\]
as $\Lambda_{\cO}$-modules. There is an automorphism of $\Lambda_{\cO}$ induced by the map $\gamma \mapsto \kappa(\gamma)^{-1} \gamma$. If $X$ is any finitely generated, torsion $\Lambda_{\cO}$-module, that automorphism send the characteristic ideal of $X$ to the characteristic ideal of $X \otimes \kappa$. It follows that we can choose $\theta_{\psi, \delta}$ so that 
\[
\varphi(\theta_{\psi, \delta}) ~=~ \varphi \kappa^{-1} (\theta_{D^*}) ~=~ \varphi^{-1}\kappa( \theta_{\chi, \ep})
\]
for all $\varphi \in \Hom(\Gamma, \overline{\QQ}_p^{\times})$. In particular,  if taking $\varphi=\kappa^s$, we have 
\begin{equation}\label{eqn:fe}
\kappa^{1-s}(\theta_{\chi, \ep}) ~=~ \kappa^s(\theta_{\psi, \delta})
\end{equation}
for all $s \in \ZZ_p$. One especially interesting case is when $d=2$. Then $d^{+}(\chi)$ and $d^{+}(\psi)$ are both equal to 1. Assume that $\epp=\omega$. 
Then $\Delta_{\gothp}$ acts on $U/U^{(\dpp)}$ trivially.  
As discussed in remark 4.3, it follows that $\varphi_0(\theta_{\psi, \delta})=0$. Therefore, taking $s=0$ in (\ref{eqn:fe}), we find that $\kappa(\theta_{\chi, \omega})=0$ too.   

\bigskip

\noindent {\bf Remark 4.9.}   We return to the case where $d^{+}=1$.  The exact sequence in proposition 3.2 can be extended to $\QQ_{\infty}$ by taking direct limits. One obtains the following exact sequence of discrete $\LO$-modules. 
\begin{equation}\label{eqn:basicexactseq}
0 \longrightarrow  H^1_{unr}(\QQ_{\infty}, D) \longrightarrow S_{\chi, \varepsilon}(\QQ_{\infty}) \longrightarrow \Hom_{\cO[\Delta]}\Big( \big(\cU_{\infty,p, \cO}^{[\varepsilon]}\big)^{(\chi)} \big/ (\overline{U}_{\infty,p,\cO}^{[\ep]})^{(\chi)}, ~ D \Big) \longrightarrow 0~~.
\end{equation}
The intervening $\LO[\Delta]$-modules in the above sequence are defined by 
\[
 \cU_{\infty, p, \cO}^{(\ep)} = \prod_{\gothp | p} ~ \cU_{\infty, \gothp, \cO}^{(\epp)}~, \quad \quad \quad \overline{U}_{\infty,p,\cO}^{[\ep]}~=~ \varprojlim_n ~ \overline{U}_{n,p,\cO}^{[\ep]}~~,
\]
where the inverse limits are defined by the norm maps. One can then take the $\chi$-components of those modules which will be finitely generated as $\LO$-modules. 
\smallskip

The term $H^1_{unr}(\QQ_{\infty}, D)$ in (\ref{eqn:basicexactseq}) can be identified with 
\[
\Hom_{\cO[\Delta]}\Big( \big(\Gal(L_{\infty}/K_{\infty})\otimes_{\ZZ_p}\cO\big)^{(\chi)},~ D\Big)~~ 
\]
where $L_{\infty}$ is the maximal abelian unramified pro-$p$-extension of $K_{\infty}$. Now $\Gal(L_{\infty}/K_{\infty})$ is a finitely generated, torsion $\Lambda$-module, a famous theorem of Iwasawa.  An immediate consequence is that 
the Pontryagin dual of  $H^1_{unr}(\QQ_{\infty}, D)$ is a finitely generated, torsion $\LO$-module.  Let $\theta_{unr, \chi}$ be a generator of its characteristic ideal. It is worth noting that the characteristic ideal for the $\LO$-module   $\big(\Gal(L_{\infty}/K_{\infty}\big)\otimes_{\ZZ_p}\cO\big)^{(\chi)}$ is then generated by $\theta_{unr, \chi}^d$.  
\smallskip

One can derive another exact sequence which is somewhat simpler than (\ref{eqn:basicexactseq}).  It will lead to the interpolation property stated in the introduction under certain assumptions.  We will discuss this carefully in the sequel to this paper and just state it here. Choose one prime $\gothp$ lying above $p$.

 We use the notation  $\big|\overline{U}_{\infty,\gothp,\cO}\big|_{\gothp}$ to be the image of $(\overline{U}_{\infty,p,\cO}^{[\ep]})^{(\chi)}$ under the projection map from $\cU^{(\ep)}_{\infty, p, \cO}$ to $\cU_{\infty, \gothp, \cO}^{(\epp)}$. Of course, it depends on $\chi$ and $\ep$ as well as $\gothp$. The exact sequence is

\begin{equation}\label{eqn:simplerbasicseq}
0 \longrightarrow  H^1_{unr}(\QQ_{\infty}, D) \longrightarrow S_{\chi, \varepsilon}(\QQ_{\infty}) \longrightarrow \Hom_{\cO}\Big( \cU_{\infty,\gothp, \cO}^{(\epp)} \big/ \cE_{\infty,\gothp,\cO}^{[\chi,\ep]}, ~ D^{(\epp)} \Big) \longrightarrow 0
\end{equation}

\noindent Note that the $\Hom_{\cO}$ term in the above exact sequence is isomorphic to the Pontryagin dual of  
$\cU_{\infty,\gothp, \cO}^{(\epp)} \big/ \cE_{\infty,\gothp,\cO}^{[\chi,\ep]}$ as a $\LO$-module. 

\bigskip

\section{The $p$-adic $L$-function of an Artin representation}

We will discuss the construction of a $p$-adic $L$-function for the representation $\rho$ in this section.  We assume that $d=2$.   The construction relies on Hida theory and on the existence of a 2-variable $p$-adic $L$-function.  As before, we choose a finite set $\Sigma$ of primes containing $p$ and $\infty$  and all the primes ramified in $K/\QQ$.  We can then regard $\rho$ as a 2-dimensional representation of $\Gal(\QQ_{\Sigma}/\QQ)$. 
The basic
idea is to embed $\rho$ into a $p$-adic family of 2-dimensional representations $\rho_k$ of $\Gal(\QQ_{\Sigma}/\QQ)$ for which a $p$-adic $L$-function
is known to exist, and then take the limit as $\rho_k$ tends towards $\rho$.  Roughly speaking, for $k \ge 2$, $\rho_k$ will be the $p$-adic representation corresponding to a modular form of weight $k$ which is ordinary at $p$ and the corresponding $p$-adic $L$-function is now classical. (See \cite{mtt}.)

Since $d^{+}=1$, it follows that $\rho$ is odd.  The Artin conjecture then implies that $\rho$ is associated to a modular form $f$ of weight 1 for $\Gamma_1(N)$, where $N$ is the Artin conductor of $\rho$. Thus, $f$ is a normalized newform of weight 1 with $q$-expansion $f=\sum a_n q^n$  and level $N$. By work of Langlands, Tunnell, Weil, and Hecke,  this is known whenever the image of $\rho$ inside $GL_2(\cO)$ is a solvable group.  In the non-solvable case (i.e., when $\rho$ is {\em icosahedral}), there are also partial results of Buzzard, Dickinson, Sheppard-Barron, and Taylor. For all of this, see \cite{bdst} and the references there. 
We assume from now on that  $\rho$ is modular of level $N$ in the above sense.

Recall that we have fixed embeddings $\sigma_{\infty}$ and $\sigma_p$ of $\overline{\QQ}$ into $\CC$ and into $\overline{\QQ}_p$, respectively.  The coefficients $a_n$  of $f$ are algebraic integers in $\CC$ and can be regarded as elements of $\overline{\QQ}_p$.    In addition, if $q$ is any prime and  $\gothq$ is a prime of $K$ lying above $q$, and $\Delta_{\gothq}$ denotes the corresponding decomposition subgroup of $\Delta$,  then the restriction of a Frobenius automorphism for $\gothq$ to the maximal subspace of $V$ on which the inertia subgroup for $\gothq$ acts trivially has trace $a_q$.   Thus, the $a_q$'s are in $\cO$ for all primes $q$. Furthermore, for $q \nmid N$, it is clear that  $\det\big(\rho(\frob_{\gothq})\big)$ is in $\cO^{\times}$.  It follows that $a_n \in \cO$ for all $n \ge 1$.

Since Hypothesis $A$ is assumed to be satisfied,  for our fixed prime $\gothp$ of $K$ lying over $p$ and our choice of $\epp$,  there is an exact sequence 
\begin{equation}\label{ord}
0\longrightarrow \cO(\varepsilon_{\gothp} )\longrightarrow T  \longrightarrow
\cO(\varepsilon'_{\gothp})\longrightarrow 0, 
\end{equation} 
of $\cO[\Delta_{\gothp}]$-modules,   where $\varepsilon_{\gothp} \neq \varepsilon'_{\gothp}$. 
As in the introduction, $\varepsilon_{\gothp}$ and  $\varepsilon'_{\gothp}$ are determined by certain characters $\varepsilon$ and
$\varepsilon'$ of $\cG = \Gal(\cK/\QQ_p)$.  
We assume at first that $\rho$ is ordinary in the sense that $\varepsilon'$ is  unramified.  
Then there exists a modular eigenform $f_1=\sum b_n q^n$ of weight 1 
with the property that $b_n = a_n$ if $(n, p)=1$, but $b_p=\varepsilon'(p)$. The form $f_1$ will be a $p$-stabilized newform in the sense of \cite{wil88}, meaning that it has level $N$ if $p\vert N$ or level $Np$ otherwise. If $p$ divides $N$, then $f=f_1$. Since the values of $\ep$, and hence $\ep'$, are in $\cO$, it follows that the $b_n \in \cO$ for all $n \ge 1$.

The maximal ideal of $\cO$ is $p\cO$. Since $\rho$ is  absolutely irreducible, and $p \nmid |\Delta|$, it follows that the residual representation $\overline{\rho}$ obtained by the canonical reduction map $GL_2(\cO)\rightarrow GL_2(\cO/p\cO)$ is   also absolutely irreducible.  Using fundamental ideas due to Hida, Wiles has constructed a deformation of $\rho$ in this setting, whose description we briefly recall.   For background in Hida's theory, we
refer the reader to \cite{hid85}, \cite{wil88}, and \cite{epw}.   

\medskip  {\bf (a) The local ring $\mathbf{h}_{\gothm}$:}   Let $\mathbf{h}_{\infty}$ denote Hida's universal ordinary Hecke algebra of level $M$, where $M$ is the prime-to-$p$ part of $N$.  It is a finite flat algebra over the subring $\Lambda' = \ZZ_p[[\Gamma']]$, where $\Gamma' \cong 1+p\ZZ_p$, considered as a subgroup of the group of diamond operators,  and so $\Lambda'$ is isomorphic to a formal power series ring  $\ZZ_p[[T]]$.
Now $f_1$ determines a maximal ideal $\gothm$ of $\mathbf{h}_{\infty}$; we write $\mathbf{h}_\gothm$ for the completion of $\mathbf{h}_\infty$ at $\gothm$.     Since the coefficients of $f_1$ generate the ring $\cO$   and that ring is generated by roots of unity of order prime to $p$, one sees that $\hm$ contains $\cO$ as a subring. Thus, $\hm$ also contains the subring $\LO' \cong \cO[[\Gamma']]$.

\medskip  {\bf (b) The Galois representation $\breve\rho~$:}  It follows from the work of Mazur, Ribet, Wiles, and others (see \cite{wil95}, Ch. 2) that the localization $\hm$ is a
Gorenstein ring. Hida's construction then implies that there exists a representation $\breve\rho: \Gal(\QSQ) \rightarrow GL_2(\hm)$.  Wiles has shown in \cite{wil88} that there exists a weight-one prime ideal $P_1=P_{\rho}\in  \text{Spec}(\mathbf{h}_{\gothm})$ such that $\breve\rho$ specializes to $\rho$ at the point $P_\rho$. If $q \nmid Np$, then the trace of a Frobenius element $\frob(q)$ for a prime above $q$  is given by the Hecke operator
$T_q$ regarded as an element of $\hm$. Furthermore, we identify $G_{\QQ_p}$ with a subgroup of $G_{\QQ}$ by the fixed embedding $\sigma_p$ (which determines $\gothp$).
Let $\mathbf{T}$ denote a
realization of $\breve\rho$. Then $\mathbf{T}$ is a free $\hm$-module of rank 2 and  there exists a $G_{\QQ_p}$-stable $\hm$-submodule
$\mathbf{T}_0\subset\mathbf{T}$ such
that $\mathbf{T/T}_0$ is unramified at $p$.   Both $\mathbf{T}_0$ and $\mathbf{T/T}_0$ are free $\hm$-modules of rank 1.  The eigenvalue of   
$\frob(p)$ on $\mathbf{T/T}_0$ is given by the Hecke operator $U_p$. Furthermore, the residual representations for $\mathbf{T}_0$ and $\mathbf{T/T}_0$ are $\epp$ and $\epp'$, respectively. 
The
determinant of $\breve\rho$ is a homomorphism $\Gal(\QSQ) \rightarrow \LO'^{\times}$.

\medskip
 {\bf (c) Specialization:} 
Recall that a point $\phi: R\rightarrow \mathbf{C}_p$ of a
$\Lambda'$-algebra $R$ is said to be arithmetic of weight $k \in \ZZ$ if the restriction $\phi |_{\Lambda'}$ is the $\cO$-algebra homomorphism 
 induced by a character $\ZZ_p^{\times} \rightarrow \overline{\QQ}^{\times}_p$ of the form 
$a\mapsto \psi(a) a^k$, where $\psi$ is a
finite-order character.  Note that if $\phi$ has values in $\cO$, then $\psi$ must have order prime to $p$. For each arithmetic point
$\phi_k\in\text{Spec}(\mathbf{h}_{\mathfrak m})$ of weight $k\geq 2$, the
specialization of $\breve\rho$ to $\phi_k$ is the representation
associated by Deligne to a certain $p$-ordinary modular form
$f_k=f_k(\phi_k)$.

\medskip  {\bf (d) $\Lambda'$-adic forms:} Finally, we recall
that a $\Lambda'$-adic form $\mathfrak{F}$ of tame level $M$ with coefficients
in a finite 
 $\Lambda'$-algebra $R$ is defined to be a formal
power series
\[
\mathfrak{F}= \sum a_n(\mathfrak{F})q^n ~ \in ~R[[q]] 
\]
 such that, for almost all
arithmetic points $\phi_k\in \text{Spec}(R)$ with $k\geq2$, the
specialization $\mathfrak{F}(\phi_k)$ is the $q$-expansion of a $p$-stabilized
newform of weight $k$ and level $Mp^r$, for suitable $r$. We say that
the form $\mathfrak{F}$ contains the classical form $f_k$ of weight $k\geq 1$
and level $Mp^r$ if there exists an arithmetic point $\phi_k$ of
weight $k$ such that $\mathfrak{F}(\phi_k)$ is the $q$-expansion of $f_k$. 
The form $\mathfrak{F}$ is said to be \emph{ordinary} if $e\mathfrak{F}=\mathfrak{F}$,
where $e$ is Hida's ordinary projector. From now on, we only consider ordinary forms. In 
the ordinary case,  a $\Lambda'$-adic eigenform with coefficients in $R$ is the same thing as a
$\Lambda'$-algebra homomorphism $\mathbf{h}_\infty\rightarrow R$.  We will assume that the homomorphism factors through $\hm$.  We can assume that $R$ is generated by the coefficients of $\mathfrak{F}$ and hence that the homomorphism is surjective.  We will also generally assume that $R$ is a domain and so $R \cong \mathbf{h}_{\infty}/\gotha$, where $\gotha$ is a prime ideal in $\mathbf{h}_{\infty}$ of height 0 (i.e., a minimal prime ideal of $\mathbf{h}_{\infty}$). Thus, $R$ is an irreducible component in the $\Lambda'$-algebra $\mathbf{h}_{\infty}$ . It obviously has characteristic 0 and contains $\Lambda'$ as a subring. 
\medskip

We write $\mathfrak{M}_R$ for the module of ordinary $\Lambda'$-adic modular forms with
coefficients in the $\Lambda'$-algebra $R$. Then $\mathfrak{M}_R$
is a finite free $R$-module, and, if $\mathfrak{M}=\mathfrak{M}_{\Lambda'}$, then there
is an isomorphism $\mathfrak{M}_R=\mathfrak{M}\otimes_{\Lambda'}R$. The ring
$\mathbf{h}_\infty$ acts on   $\mathfrak{M}$, and the module $(\mathfrak{M}_R)_\gothm$ is
canonically isomorphic to the $\mathbf{h}_\gothm$-module
$\Hom_R(\mathbf{h}_\gothm\otimes_{\Lambda'}R, R)$.
\medskip

With these notions in hand, we can state a key result.   

\bigskip
 
\noindent {\bf Theorem 5.1.} (Wiles)  {\em Suppose that $f_1$ is the $p$-stabilized form associated to $\rho$ and the unramified character $\ep'$.  Then there exists an ordinary $\Lambda'$-adic 
eigenform $\mathfrak{F}$, defined
over an extension $R$ of $\Lambda'$, which specializes to $f_1$ in weight 1.} 

\bigskip

Note that in the above theorem we can assume that $R$ is a domain and is an irreducible component in the $\LO'$-algebra $\hm$. Now $\ep \neq \ep'$ since we are assuming that Hypothesis A holds for $\rho$. 
If $\ep$ is unramified, then one can reverse the roles of $\ep$ and $\ep'$, obtaining a different $\Lambda'$-adic 
eigenform (which would correspond to a different choice of $\gothm$ and a different $R$).   However, for our fixed choice of $\ep$, there is a uniqueness result which is very important for our purpose.  
\bigskip

\noindent {\bf Theorem 5.2.} (Greenberg-Vatsal, Bell\"aiche-Dimitrov)  {\em Suppose that a choice 
of $\varepsilon$ is fixed and that Hypothesis A holds. Then the form $\mathfrak{F}$ is unique.}

\bigskip

\noindent In particular, $\gothm$ and $R$ are uniquely determined by $\chi$ and $\ep$.    A version of this theorem was discovered by the authors of this paper in 1997, but never published. We needed the restrictive assumption that $\ep'/\ep$ is not of order 2. 
It was subsequently rediscovered by Bell\"aiche-Dimitrov \cite{bd16}, who actually proved the result without that assumption. 
\medskip

It is important to observe that the residue 
ring of $\mathbf{h}_\gothm$ at the prime
ideal $P_\rho$ is exactly equal to $\cO$, rather than some larger ring $\cO'$.  Recall that $\cO$ is generated over $\ZZ_p$ by the values of $\chi$ and $\ep$, the ring of integers in $\cF=\QQ_p(\chi, \ep)$.   We justify the claim as follows. The Hecke algebra $\mathbf{h}_\gothm$ is generated by the Hecke operators
and the diamond operators. We must show that the images of these generators are in $\cO$. 

For the diamond operators, this is clear since the values of $det \circ \rho$ are in that ring. 
Now the image of any Hecke operator modulo $P_\rho$ is the eigenvalue
of the corresponding Hecke operator for the weight 1 form associated to $\rho$. We have
to show that these eigenvalues lie in $\cO$. For any prime $q$, let $\gothq$ be a prime of $K$ lying over $q$. 
 For the Hecke operator
$T_q$ with $q\nmid Np$, the eigenvalue of $T_q$ is the trace of $\rho(\frob(\gothq))$, which
lies in $\cO$ by definition.    If $q$ divides $Np$, then we argue as follows.  As explained in the introduction, $\rho$ can be realized in a 2-dimensional vector space $V$ over the field $\cF$.  If $q \nmid p$, let $I_{\gothq}$ be the inertia subgroup of $\Delta$ for the prime $\gothq$.  
Then $V^{I_\gothq}$ is an $\cF$-subspace of $V$ of dimension $\le 1$. If $V^{I_\gothq}$ is 1-dimensional,
then the eigenvalue of $U_q$ is equal to the eigenvalue of Frobenius on that subspace. That eigenvalue obviously lies in $\cF$, and hence in $\cO$.  If that subspace has dimension 0, then the image of $U_q$ is 0. Finally, $U_p$ acts on $f_1$ with eigenvalue $\epp'(\frob(\gothp))$. Since $\chi |_{\DP}= \epp + \epp'$, it follows that the values of $\epp'$ are in $\cO$. Thus, all the eigenvalues of the Hecke operators are indeed in $\cO$. 
\medskip

\bigskip
The above discussion yields a representation $\breve\rho:  \Gal(\QSQ) \rightarrow GL_2(\mathbf{h}_\gothm)$
which specializes to $\rho$ at an arithmetic specialization of weight 1. We can compose with the surjective homomorphism $\hm \rightarrow R$ to obtain a representation 
\[
\tilde\rho:~ \Gal(\QSQ)~ \longrightarrow ~ GL_2(R)~~.
\]
Furthermore, there is a homomorphism $\phi_{\rho}: R \to \cO$ such that $\phi_{\rho}\circ\tilde\rho=\rho$.  Here we are also using the notation $\phi_{\rho}$ for the induced continuous  group homomorphism 
$GL_d(R) \rightarrow GL_d(\cO)$.  Such representation $\breve\rho$ and $\tilde\rho$ are called deformation of $\rho$. The rings $\hm$ and  $R$ are the corresponding deformation rings. Each deformation has infinitely many motivic points, namely the points corresponding to forms of higher weight.  The deformation $\tilde\rho$ is the one we will concentrate on. 
These motivic points
have well-defined $p$-adic $L$-functions (up to multiplying by units in the corresponding residue rings of $R$). We can therefore define a $p$-adic $L$-function
for $\rho$ by continuity.  We now explain this more precisely. 
\medskip

Thus we are led to describe the $p$-adic $L$-function associated to the
representation $\tilde\rho$. It will correspond to an element in $R[[\Gamma]]$,  where $\Gamma=\Gal(\QQ_{\infty}/\QQ)$ is as in the previous sections. This $L$-function was constructed by
Kitagawa \cite{kit94} and by Greenberg-Stevens
\cite{gs}.  We want to recall the modular-symbol construction of
Greenberg-Stevens, as this is the most convenient for us. Our discussion here is only a summary, and we refer the
reader to the original paper of Greenberg and Stevens for a detailed
exposition. For each positive integer $r$, let $\Gamma_r$ denote
the group $\Gamma_1(Np^r)$. Let $\mathbf
M$ denote the inverse limit of the cohomology groups
$e\hp(\Gamma_r, \ZZ_p)$; then $\mathbf M$
is a module over the Hecke algebra $\mathbf{h}_\infty$. The matrix
$\iota=\left(\begin{array}{cc}1&0\\ 0&-1\end{array}\right)$ induces a
Hecke-equivariant involution on $\mathbf M$, and we write ${\mathbf
M}^{\pm}$ for the $\pm$-eigenspaces under this action. Evidently we
have the decomposition ${\mathbf M}={\mathbf M}^+\oplus{\mathbf M}^-$. Since the maximal
ideal $\gothm$ corresponds to an irreducible representation, it follows
from Thm. 4.3 and Lemma 6.9 of \cite{gs} that we have an
isomorphism of $\mathbf{h}_\gothm$-modules
$${\mathbf M}_{\gothm}=\text{Symb}_{\Gamma_1}({\mathbf D})_{\gothm},$$  where
$\mathbf D$ is the space of $\ZZ_p$-valued measures on the set of
primitive elements of $\ZZ_p^2$, and $\text{Symb}_{\Gamma_1}({\mathbf
D})$ denotes the group of modular symbols over the group $\Gamma_1$
with values in $\mathbf D$. Furthermore, since $\mathbf{h}_\gothm$ is
Gorenstein, a duality argument shows that, for any choice $\alpha=\pm$
of the sign, there are isomorphisms of Hecke modules
$\mathfrak{M}_{\gothm}={\mathbf{M}}_{\gothm}^{\alpha}=\mathbf{h}_\gothm$. 

Now let $\mathfrak{F}$ be a $\Lambda'$-adic eigenform form, defined over the ring $R$,
such that $\mathfrak{F}$ specializes to $f_1=f_{\rho}$ in weight 1. 
The choice of $\epp'$ made in (\ref{ord}) is implicit in the choice of $\mathfrak{F}$. 
Thus we
have a homomorphism $\mathbf{h}_\gothm\rightarrow R$ such that the kernel
$\mathfrak{P}$ is a minimal prime contained in the height 1 prime ideal 
 $P_\rho$ of $\hm$.
Extending scalars via
$\Lambda'\rightarrow R$, we obtain isomorphisms 
$\mathfrak{M}_\gothm\otimes_{\Lambda'} R=\mathbf{M}^\pm_\gothm\otimes_{\Lambda'} R=\mathbf{h}_\gothm\otimes_{\Lambda'} R$. 
Let $\delta^{\alpha}(\mathfrak{F})$ be
the image of $\mathfrak{F}\in\mathfrak{M}\otimes_{\Lambda'} R$ in the module
${\mathbf{M}}_\gothm^\pm\otimes_{\Lambda'} R$, and let $\delta(\mathfrak{F})=\delta^{+}(\mathfrak{F})+
\delta^{-}(\mathfrak{F})$. Finally, let 
$\mathfrak{L}(\mathfrak{F})=\delta(\mathfrak{F})(\{0\}-\{\infty\})\in\mathbf{D}\otimes_{\Lambda' }R$ denote the
special value (\cite{gs} Sec. 4.10). 
Then $\mathfrak{L}(\mathfrak{F})$ is an element of $\mathbf{D}
\otimes_{\Lambda'}R$, and the construction of \cite{gs}, equations 5.2a and 5.4,
gives an element
\begin{equation}
\label{2var} 
L_p(\mathfrak{F})\in R[[\mathbf{Z}_p^\times]]\end{equation} 
known as the standard two-variable $p$-adic $L$-function
associated to $\mathfrak{F}$.  The terminology is justified as follows. 
Let $P_k$ denote an arithmetic point in
$\text{Spec}(R)$, of weight $k\geq 2$. Then $P_k$ induces a specialization
map $\phi_k:R\rightarrow\cO$,
as explained above.  We may specialize  $L_p(\mathfrak{F})$
via $\phi_k$ to obtain a measure $L_p(\mathfrak{F})(P_k)$ on $\mathbf{Z}_p^\times$; one
sees from the definitions that, if $f_k$ denotes the specialization of
$\cF$ to $P_k$, then the measure $L_p(\mathfrak{F})(P_k)$ is a $p$-adic $L$-function
for $f_k$. Specifically, if $\psi$ is any finite order character of $\mathbf{Z}_p^\times$,
then
$\mu_k=L_p(\mathfrak{F})(P_k)$ satisfies the property
\begin{equation}\label{eqn:smooth}
\int_{\mathbf{Z}_p^\times} \psi d\mu_k=e_p(f_k, \psi) \tau(\psi)
\frac{L(f_k\otimes\psi^{-1}, 1)}{(-2\pi i)\Omega_{k}^\alpha}
\end{equation}
for certain periods $\Omega^{\pm}_k$.
Here  $\tau(\psi)$ is the standard Gauss sum for $\psi$
and 
$$e_p(f_k, \psi)= \alpha_k^{-n}(1 - \alpha_k^{-1}\psi(p)),$$
where $\alpha_k$ is the (unit) eigenvalue of the $U_p$ operator on $f_k$. The
quantity $L(f_k\otimes\psi^{-1}, 1)$ denotes the standard $L$-function of $f_k$, 
twisted by the finite order character $\psi^{-1}$ and evaluated at $s=1$. 
A similar but more complicated formula holds for the integral of $\psi z^{j-1}$, where 
$j$ is any integer in the critical range $1\leq j\leq k-1$. It is essentially the formula in the proposition of section 14 in \cite{mtt}.

It follows easily
from Hida's control theorems \cite{hid85} 
that the periods $\Omega_k^{\alpha}$
are the {\it canonical periods} of $f_k$ introduced in \cite{vat97}. These periods
are defined up to $p$-adic units. The sign 
$\alpha=\pm$ in the period is determined by the parity of $\psi$. 

We can identify $\ZZ_p^{\times}$ with $\Gal(\QQ(\mu_{p^{\infty}})/\QQ)$ in the usual way.  Then $1+p\ZZ_p$ is identified with $\Gamma$. The characters $\psi$ of $\ZZ_p^{\times}$ can be regarded as characters of that Galois group.  We have $\psi =\omega^t \varphi$ where $0 \le t \le p-2$ and $\varphi$ is a character of finite order of $\Gamma$. 

The ring $R[[\ZZ_p^{\times}]]$ has a natural $R$-algebra involution $\iota$ defined by $\iota(z)=z^{-1}$ for all $z \in \ZZ_p^{\times}$. For our purpose, we consider the element $\iota\big(L_p(\mathfrak{F})\big)$ in $R[[\ZZ_p^{\times}]]$. If we specialize via $\phi_k$ as above, we obtain a measure $\mu_k^{\iota}$. Integrating $\psi$ against   $\mu_k^{\iota}$ gives a formula like that in (\ref{eqn:smooth}), but with $\psi$ replaced by $\psi^{-1}$ on the right side. In particularly, the $L$-value occurring in the formula would now be $L(f_k \otimes \psi, 1) =L(f_k, \psi,1)$.  Now $\ZZ_p^{\times} \cong \FF_p^{\times} \times (1+p\ZZ_p)$ and so one has the decomposition 
\[
R[[\ZZ_p^{\times}]] ~ \cong ~ \bigoplus_{t=0}^{p-2} ~ e_{\omega^t}R[[\ZZ_p^{\times}]]
\]
where $e_{\omega^t}$ is the idempotent for $\omega^t$.  Each summand is isomorphic to $R[[1+p\ZZ_p]]$ which we identify with $R[[\Gamma]]$.  If we project $\iota\big(L_p(\mathfrak{F})\big)$ to the $\omega^t$-component for a fixed $t$, we get an element $\Theta_t \in R[[\Gamma]]$. Specializing via $\phi_k$, we get an element $\theta_{t,k} \in (R/P_k)[[\Gamma]]$. 
If we let $\psi = \omega^t \varphi$, where $\varphi$ is a character of $\Gamma$, then integrating $\psi$ against $\mu_k^{\iota}$ gives the value $\varphi(\theta_{t,k})$. Note that $L(f_k \otimes \psi, 1)=L(f_k\otimes \omega^t, \varphi, 1)$. 

Recall that we assumed earlier that $\rho$ is ordinary in the sense that there is an exact sequence (\ref{ord}) with $\epp'$ unramified. However, in general, the action of $G_{\QQ_p}$ on $T/\cO(\epp)$ might be ramified. One could then replace $\rho$ by $\rho \otimes \omega^{-t}$ for some $t$, $0 \le t \le p-2$, to obtain an ordinary Artin representation.  One can apply Wiles theorem to the corresponding weight 1 form $f_1$ and then tensor by $\omega^t$ obtaining the family $f_k \otimes \omega^t$ of modular forms which specializes in weight 1 to $f_{\rho} = f_1\otimes \omega^t$. With this motivation in mind, it is natural to define $\theta^{an}_{\chi, \ep}$ to be the image of $\Theta_t$ by the specialization map $\phi_{\rho}:  R[[\Gamma]] \to \cO[[\Gamma]] = \LO$. This is the element alluded to in the introduction and is well-defined up to multiplication by an element of $\cO^{\times}$.

\bigskip

\noindent {\bf Remark 5.3.} 
As we have remarked several times already, the $p$-adic $L$-function 
$L_p(s, \chi, \varepsilon)$
depends on the choice of $\varepsilon$, in the case $(N, p)=1$. 
Furthermore, the $L$-function is only defined up to a $p$-adic unit, owing to 
the indeterminacy of the complex periods in the definition of the 2-variable
$p$-adic $L$-function. Finally,  it is {\it not at all} clear that $L_p(s, \chi, \varepsilon)$ is nonzero. 
However, in view of Theorem 3 of the introduction, it seems reasonable to expect
that $L_p(s, \chi, \varepsilon)$ should be related to the Coleman series of the norm-compatible 
family of units given in that theorem, and that a formula in terms of logarithms of global units similar to the one in Theorem 3 should hold. 
We will take up this subject in the sequel to this paper. 
\bigskip

\noindent {\bf Remark 5.4.} The modular cuspforms of weight one may be classified according to the isomorphism class
of the image of the associated representation $\rho$ in $PGL_2(\mathbf{C})$. The possible image $G$ of $\Delta$
is isomorphic to one of the following: a dihedral group, $A_4$, $S_4$, or $A_5$. In the event that $G$
is dihedral, it is induced from a character of a quadratic extension of $\mathbf{Q}$. Then $\rho$ is called
CM (complex multiplication) or RM (real multiplication), according to whether it is induced from an imaginary quadratic field or not. An
alternative approach to the $p$-adic $L$-function in this case is via the Katz 2-variable $p$-adic $L$-function; the connection
with global elliptic units is then given by the $p$-adic Kronecker limit formula (see \cite{deshalit}, Theorem 5.2 and the subsequent discussion). 
We refer the reader to Ferrara \cite{ferrara-thesis} for more details. When $\rho$ has real multiplication, Ferrara has also compiled numerical evidence
relating $\phi(\theta^{an}_{\chi,\ep})$ to global units (in fact, Stark units) in abelian extensions of real quadratic fields when $\phi$ has order $p$.   On the other hand, when $\rho$ is an exotic form of 
non-dihedral type, absolutely nothing is known. 
\bigskip

\noindent {\bf Remark 5.5.} As stated in the introduction, the natural formulation of an Iwasawa main conjecture is the assertion that  $\theta_{\chi,\ep}^{al}$ and $\theta_{\chi,\ep}^{an}$ generate the same ideal in $\LO$.  As mentioned in remark 5.3, we cannot even rule out the possibility that $\theta_{\chi,\ep}^{an}=0$.  The fact that $\theta_{\chi,\ep}^{al} \neq 0$ is the content of proposition 4.2.  It seems reasonable to believe that the main conjecture in weight 1 should follow  from the main
conjecture in higher weights.  Some ideas in \cite{oc} would possibly be useful in showing that. However, regrettably,
the proof of the main conjecture in higher weight
by Urban-Skinner \cite{urban-skinner} does not apply
to the case at hand because those authors require 
the presence of a prime $q\vert M$ such that $q$ is a prime of multiplicative reduction for $\mathfrak{F}$. Theorem 5.6 below implies that the form $\mathfrak{F}$ admits no such prime $q$.  On the other hand, it appears that the existence of such a $q$ is only
required to ensure that a certain anticyclotomic $\mu$-invariant vanishes, and this
can also be proven when $q$ is a prime of supercuspidal reduction, so it may be that
the result of Skinner-Urban applies when there exists a prime $q$ such that the image
of a decomposition group at $q$ under $\rho$ is non-abelian and irreducible. We have 
not pursued this idea.
\bigskip

We have assumed $d=2$ until now. We close this paper with some general comments for arbitrary $d$. Let $R$ be a complete Noetherian local ring with finite residue field of characteristic $p$.   Assume that we have 
a continuous representation  $\tilde\rho: \Gal(\QQ_{\Sigma}/\QQ) \rightarrow GL_d(R)$, together with a
 continuous, surjective ring homomorphism $\phi_{\rho}: R \rightarrow \cO$ such that
$\phi_{\rho}\circ\tilde\rho=\rho$.  That is, $\tilde\rho$ is a deformation of $\rho$. The ring $R$ is the corresponding deformation ring. Since $\cO$ is generated over $\ZZ_p$ by roots of unity of order prime to $p$, it is easy to see that $R$ contains a subring which we can identify with $\cO$. Therefore, $R$ is an $\cO$-algebra and $\phi_{\rho}$ is an $\cO$-algebra homomorphism. Furthermore, the residue ring of $R$ is the same as that of $\cO$, namely $\cO/p\cO$.

 We write $P_{\rho}$ for the kernel
of the ring homomorphism $\phi_{\rho}$, a prime ideal in $R$.  If $P$ is any other prime ideal of $R$ such that $R/P \cong \cO'$, where $\cO'$ is a finite integral extension of $\cO$, then composing $\tilde\rho$ with the group homomorphism $GL_d(R) \to GL_d(\cO')$ defined by reducing modulo $P$ gives a $d$-dimensional representation of $ \Gal(\QQ_{\Sigma}/\QQ)$ over $\cO'$. Thus, we get a family of representations indexed by such prime ideals $P$ of $R$. 
Of course, $\tilde\rho$ mod $P$ is not necessarily an Artin representation. However, if $R$ is a domain, then  the following proposition implies that it has one property in common with $\rho$. Namely,  $\tilde\rho$ mod $P$ is potentially unramified at all primes not dividing $p$. A similar result in the case $d=2$ is contained in theorem 2.3 in \cite{oc}.
 \bigskip
 
\noindent {\bf Theorem 5.6.}  {\em  Let $\tilde\rho$ be a deformation of the Artin representation $\rho$. Assume that
the corresponding deformation ring $R$ is an integral domain. Suppose that $q$ is a finite prime in $\Sigma$ and that $q \neq p$. Let $I_v$ be the inertia subgroup of $\Gal(\QQ_{\Sigma}/\QQ)$ for a prime $v$ of $\QQ_{\Sigma}$ lying over $q$. Then the image of $I_v$ under $\tilde\rho$ is finite. }
\bigskip

\noindent {\em Proof.}    Recall that $\rho$ factors through $\Gal(K/\QQ)$. Thus, the residual representation for $\tilde\rho$ also factors through $\Gal(K/\QQ)$. It follows that the image of $\Gal(\QQ_{\Sigma}/K)$ under $\tilde\rho$ is a pro-$p$ group and hence that $\tilde\rho$ factors through $\Gal(M/\QQ)$, where $M$ is the maximal pro-$p$-extension of $K$ contained in $\QQ_{\Sigma}$. Therefore, since $q \neq p$, the restriction $\tilde\rho|_{G_{K_v}}$ factors through $\Gal(K_v^{tr,p}/K_v)$,  where  $K_v^{tr,p}$ denotes the maximal
pro-$p$ tamely ramified extension of $K_v$.  
Enlarging $K$ if necessary, we can assume that $K$, and hence $K_v$, contains $\mu_p$.   Now $K_v^{tr,p}$ contains the
field $K_v^{unr,p}$, the maximal unramified pro-$p$ extension of $K_v$. In fact, $K_v^{unr,p} = K_v(\mu_{p^{\infty}})$ and $K_v^{tr,p} = K_v(\mu_{p^{\infty}}, \sqrt[p^{\infty}]{q})$. Furthermore, both $\Gal(K_v^{^{unr,p}}/K_v)$  and. $\Gal(K_v^{tr,p}/K_v^{unr,p})$ are isomorphic to  $\ZZ_p$.    The Frobenius element of $\Gal(K_v^{unr,p}/K_v)$  acts on $\Gal(K_v^{tr,p}/K_v^{unr,p})$ by conjugation. This action  is given by $x \mapsto x^{\cN(v)}$ for all $x \in  \Gal(K_v^{tr,p}/K_v^{unr,p})$. Here $\cN(v) = q^f$, where $f$ is the residue field degree for $K_v/\QQ_q$.   Thus,  $\Gal(K_v^{tr,p}/K_v)$ is an extension of $\ZZ_p$
by $\ZZ_p(1)$. 

Let $\tau$ denote a topological generator of $\Gal(K_v^{tr,p}/K_v^{unr,p})$. The
restriction of $\tilde\rho$ to the inertia subgroup $\Gal(K^{tr,p}_v/K_v^{unr,p})$ is determined by the matrix  $\tilde\rho(\tau)$. The eigenvalues of this matrix are in some extension of  the fraction field of $R$. In addition, if $\alpha$ is one of those eigenvalues, then so is $\alpha^{q^f}$.  Since the matrix has only finitely many eigenvalues, it follows that these eigenvalues are  roots of unity.  Therefore, for some $t \ge 1$, the matrix $A= \tilde\rho(\tau^t)$ is unipotent.  If we again enlarge the field $K$, we can assume that $t=1$ i.e., that the eigenvalues of    $\tilde\rho(\tau)$ are all equal to 1.  Of course, $f$ may change too.j

Let $\mathbf{F}$ denote the fraction field of $R$.  Thus, $\Gal(\QQ_{\Sigma}/\QQ)$ acts on $\mathbf{V}=\mathbf{F}^d$ via $\tilde\rho$.   Let $N= A-I_d$, a nilpotent matrix with entries in $R$.  
There is a filtration on $\mathbf{V}$ defined by the subspaces 
$\mathbf{W}_k=\Ker( N^k)$ for $k \ge 0$.  Here $\mathbf{W}_0 = 0$ and $\mathbf{W}_t = \mathbf{V}$ if $t$ is sufficiently large. These subspaces are $G_{K_v}$-invariant.  The action of $\Gal(K_v^{tr,p}/K_v)$ on the subquotients $\mathbf{W}_{k+1}\big/\mathbf{W}_k$ 
factors through  $\Gal(K_v^{unr,p}/K_v)$.  We will prove that $\mathbf{W}_1 = \mathbf{V}$, i.e., that $N$ annihilates $\mathbf{V}$.  

Assume to the contrary that $\mathbf{W}_1$ is a proper subspace of $\mathbf{W}_2$.  Then, 
$N\mathbf{W}_2$ is a nontrivial subspace of $\mathbf{W}_1$.   We can lift the Frobenius element of 
$\Gal(K_v^{unr,p}/K_v)$ to an element of $\Gal(K_v^{tr,p}/K_v)$ which then acts on $\mathbf{V}$ by a matrix 
$B$ with entries in $\mathbf{h}$.  
Furthermore, we have $BAB^{-1} = A^{q^f}$.  It follows that
\[
BNB^{-1}  =  q^f N +  CN^2
\]
where $C$ is a matrix with entries in $R$.  One sees easily that  $\mathbf{W}_1$  and $\mathbf{W}_2$ are invariant under multiplication by $B$ and by $N$. 
Denoting the restrictions of $B$ and $N$ to $\mathbf{W}_2$ by $B_2$ and $N_2$, respectively,  we have $B_2N_2 = q^fN_2B_2$.  One sees from this that if $\beta$ is an eigenvalue for the action of $B$ on $\mathbf{W}_2/\mathbf{W}_1$, then $q^f\beta$ is an eigenvalue for the action of $B$ on $\mathbf{W}_1$.  Thus, $B$ has two eigenvalues whose ratio is $q^f$.  

Applying the homomorphism $\phi_{\rho}$ to $B$, we obtain a matrix with entries in  $\cO$ which again has two 
eigenvalues whose ratio is $q^f$. However, this matrix is in the image of $\rho$, a finite subgroup 
of $GL_d(\cO)$, and its eigenvalues are roots of unity. This is a contradiction and so we have 
$\mathbf{W}_1 = \mathbf{V}$.  It follows that $\tau$ acts trivially on $\mathbf{V}$.  
Therefore, the image of the inertia subgroup of $G_{K}$ for $v$ under $\tilde\rho$ is trivial.  
In the proof, we have possibly replaced $K$ by a finite extension. Hence, for the original $K$,  
the image of $I_v$ is finite.     \hfill$\blacksquare$
\bigskip

In general, if we try to imitate the $d=2$ case, there are some natural restrictions on the kinds of deformations of $\rho$ we consider. We are assuming that $\rho$ satisfies Hypothesis A.  The  underlying $\cF$-representation space $V$ contains a $G_{\QQ}$-invariant $\cO$-lattice $T$.  Fix a prime $\gothp$ of $K$ lying over $p$.    Then $T^{(\varepsilon_{\gothp})} \cong \cO(\varepsilon_{\gothp})$ as $\cO[\Delta_{\gothp}]$-modules.   Since $p \nmid |\Delta_{\gothp}|$,  $T/T^{(\varepsilon_{\gothp})}$ is a free $\cO$-module of rank $d-1$.

It is natural to consider deformations of $T$ which simultaneously deform $T^{(\epp)}$.  To be precise, let $\gothp$ be the prime of $K$ induced by the fixed embedding of $\overline{\QQ}$ 
into $\overline{\QQ}_p$, as in the introduction.   We can then identify $G_{\QQ_p}$ with a subgroup of $G_{\QQ}$ and restrict $\tilde\rho$ to that subgroup.   If $\mathbf{T}$ is the underlying free $R$-module of rank 
$d$ for $\tilde\rho$, we make the following hypothesis. 
\bigskip

\noindent {\bf Hypothesis $\tilde{A}$:} 
   {\em There exists a 
$G_{\QQ_p}$-invariant $R$-submodule $\mathbf{T}_0$ of $\mathbf{T}$ with the following properties: 
\medskip

\noindent $(i)$:  $\mathbf{T}_0$ is a free $R$-module of rank 1,
\medskip

\noindent $(ii)$:   $\mathbf{T}/\mathbf{T}_0$ is also a free $R$-module, 
\medskip

\noindent $(iii)$:  the image of $\mathbf{T}_0$ under the map $\phi_{\rho}$ is 
$T^{(\varepsilon_{\gothp})}$.}
\bigskip

\noindent Note that the residual representations of $G_{\QQ_p}$ for $\mathbf{T}_0$ and for  $\mathbf{T}/\mathbf{T}_0$ are the same as for $T^{(\varepsilon_{\gothp})}$  and for $T/T^{(\varepsilon_{\gothp})}$, respectively, and that they have no irreducible constituents in common. 
\bigskip

\noindent {\bf Remark 5.7.}   Assuming that $R$ is a domain, hypothesis $\tilde\A$ is not as stringent as it might seem at first. 
As mentioned in the proof of proposition 5.6,  the image of $\Gal(\QQ_{\Sigma}/K)$ under $\tilde\rho$ is a pro-$p$ group. It follows that the action of $\Gal(\QQ_{\Sigma}/\QQ)$ on   
$\mathbf{T}$ factors through a group $\cG$ which has a normal pro-$p$ subgroup $\cN$ and corresponding quotient group $\cG/\cN \cong \Delta$. Since $|\Delta|$ is prime to $p$, this group extension splits and $\cG$ has a subgroup $\cD$ which is mapped isomorphically to $\Delta$ under the map $\cG \to \cG/\cN$.  The subgroup $\cD$ is not unique, but, for simplicity of notation, we make a choice and identify $\cD$ with $\Delta$ by that isomorphism. Then $\tilde\rho|_{\Delta}$ is closely related to  $\rho$. To be precise,  the reduction of $\tilde\rho|_{\Delta}$ modulo $P_{\rho}$ coincides with $\rho$. 
We let $\mathbf{F}$ be the fraction field of $R$ (as in the proof of proposition 5.6) and let $\mathbf{V}=\bT \otimes_{R} \mathbf{F}$. Since $R$ contains $\cO$, $\bF$ must contain $\cF$. 
The action of  $\Delta$ on the $\bF$-vector space $\bV$ is a $d$-dimensional representation over $\bF$ and the character of that representation must be $\chi$.  Now $\rho$ is a representation of $\Delta$ over $\cF$ and this $\bF$-representation of $\Delta$ must be just the extension of scalars of $\rho$ from $\cF$ to $\bF$. It then follows that, for the restriction to the subgroup $\Delta_{\gothp}$, there is a 1-dimensional $\Delta_{\gothp}$-invariant $\bF$-subspace $\bV^{(\varepsilon_{\gothp})}$ on which $\Delta_{\gothp}$ acts by $\varepsilon_{\gothp}$.  Of course, that subspace is simply $e_{\varepsilon_{\gothp}}\bV$,  where $e_{\varepsilon_{\gothp}} \in \cO[\Delta_{\gothp}]$ is the idempotent corresponding to $\varepsilon_{\gothp}$. If Hypothesis $\tilde{A}$ is satisfied, then it is clear that $\bT_0 = e_{\varepsilon_{\gothp}}\bT$.  

These observations lead to the following conclusion. It is sufficient to just assume that $\bT$ contains a $G_{\QQ_p}$-invariant $R$-submodule $\bT_0'$ of rank 1 such that the action of $\Delta_{\gothp}$ on $\bT_0'$ is given by $\varepsilon_{\gothp}$. It then follows that $\bV^{(\varepsilon_{\gothp})}$ is a $G_{\QQ_p}$-invariant $\bF$-subspace of $\bV$. 
It is then clear that one can take 
\[
\bT_0~=~ \bT \cap \bV^{(\varepsilon_{\gothp})} ~=~  e_{\varepsilon_{\gothp}}\bT~~.
\] 
It is $G_{\QQ_p}$-invariant and a direct summand in $\bT$ and  therefore is a free $R$-module.  The quotient is also free.   One has $\bT_0' \subseteq \bT_0 \subset \bV^{(\varepsilon_{\gothp})}$. Hence $\bT_0$ has rank 1. Furthermore, the image of $\bT_0$ under $\phi_{\rho}$ is just $e_{\varepsilon_{\gothp}}T$ which is $T^{(\varepsilon_{\gothp})}$.  
\smallskip

Assuming that such a $G_{\QQ_p}$-invariant $R$-submodule does exist, note that the action of $\G_{\QQ_p}$ on $\bT_0$ is given by a continuous homomorphism
\[
\tilde\varepsilon_{\gothp}: ~ G_{\QQ_p} \longrightarrow R^{\times}~~.
\]
with the property that $\phi_{\rho} \circ  \tilde\varepsilon_{\gothp} = \varepsilon_{\gothp}$. In other words, $\tilde\varepsilon_{\gothp}$ is a deformation of the $G_{\QQ_p}$-representation $\varepsilon_{\gothp}$.  
It is natural to denote $\mathbf{T}_0$ by $\mathbf{T}^{(\tilde\varepsilon_{\gothp})}$.  
It is indeed the maximal $R$-submodule of $\mathbf{T}$ 
on which $G_{\QQ_p}$ acts by $\tilde\varepsilon_{\gothp}$.  
\smallskip

In addition to Hypothesis $\tilde{A}$,  we would need to assume that the deformation $\tilde{\rho}$ of $\rho$ has a rich supply of {\em ``motivic''}  specializations so that one can reasonably expect the existence of a suitable $p$-adic $L$-function. We have no way to ensure this in general. However, the remark below describes one case for $d=3$ where we do have an ample set of motivic specializations. 

\bigskip

\noindent {\bf Remark 5.8.}  Let $\rho$ be the 
2-dimensional representation associated to a modular form of weight 1. Assume that $\rho$ factors through $\Delta=\Gal(K/\QQ)$ and that $p \nmid |\Delta|$. 
Then $\text{Ad}^0(\rho)=\text{Sym}^2(\rho)\otimes\det(\rho)^{-1}$ is a 3-dimensional representation factoring through $\Delta$
which satisfies $d^+=1$. It may or may not satisfy Hypothesis A; it will
do so if and only if $\varepsilon/\varepsilon'$ has order $\ge 3$.  If it does satisfy
Hypothesis A, then one can take the adjoint of the $\Lambda'$-adic Galois representation
associated to $\mathfrak{F}$ to get a good deformation of $\text{Ad}^0(\rho)$. The algebraic part
of our theory goes through unchanged. The analytic side also goes through,
but there is one important difference: the $p$-adic $L$-function has to be defined
via continuity from the Schmidt $L$-function constructed in \cite{coates-schmidt}.
The details require more 
care, since in the dihedral cases the symmetric square is reducible, and in the RM 
case the corresponding deformation space is not smooth over
the weight space at the weight one point of interest. However, it seems reasonable
to expect a relationship with units to hold in the CM and exotic $A_4$, $S_4$ and $A_5$
cases. We omit the details here. Regrettably, this example, and the $d=2$ example of modular
forms, are the only examples of which we are aware where $d^+=1$ and sufficiently
many motivic deformations
exist.

\end{document}